\documentclass[12pt, oneside]{amsart}   	% use "amsart" instead of "article" for AMSLaTeX format
\usepackage{geometry}                		% See geometry.pdf to learn the layout options. There are lots.
\usepackage{graphicx}				% Use pdf, png, jpg, or eps§ with pdflatex; use eps in DVI mode
\usepackage{amssymb}
\usepackage{amsthm}
\usepackage{amsmath, amssymb}

\theoremstyle{plain}
\newtheorem{Thm}{Theorem}

\newtheorem{Coro}[Thm]{Corollary}
\newtheorem{Lem}[Thm]{Lemma}
\newtheorem{Claim}[Thm]{Claim}
\newtheorem{Prop}[Thm]{Proposition}
\newtheorem{Que}[Thm]{Question}

\newtheorem{Def}[Thm]{Definition}

\begin{document}
\begin{abstract}     
We define a new integer invariant of a finite graph $\Gamma$, the {\em freeness index}, that measures the extent to which $\Gamma$ can be embedded in the 3-sphere so that it and its subgraphs have ``simple" complements, i.e., complements which are homeomorphic to a connect-sum of handlebodies.  We relate the freeness index to questions of embedding graphs into surfaces, in particular to the orientable cycle double cover conjecture.    We show that a cubic graph satisfying the orientable double cycle cover conjecture has freeness index at least two.
 \end{abstract}
\title{The Freeness Index of a Graph}
\author{Abigail Thompson}
\date{}							% Activate to display a given date or no date
 \footnote{Supported by the US-Israel Binational Science Foundation grant 2018313. The author is grateful to St. Catherine's College, Oxford, and IST Austria for their hospitality while this work was completed. }

\maketitle
\section{Introduction}\label{intro} 

A knot in the 3-sphere is the trivial knot if and only if its complement is a solid torus, or genus 1 handlebody.    We can extend the notion of triviality from knots to embedded (finite) graphs by considering a graph embedded in the 3-sphere to be ``unknotted'' if its complement is a connect-sum of handlebodies.    While every connected graph has an embedding into $S^3$ with handlebody complement, some graphs (for example $K_6$) have knotted subgraphs no matter how they are embedded into $S^3$ (\cite{CG}, \cite{S}).   Planar graphs, by contrast, can be embedded in $S^3$ such that every subgraph is unknotted.      In this paper we define an integer invariant of a finite (connected) graph $\Gamma$, the {\em freeness index}, which measures the extent to which $\Gamma$ can be embedded in the 3-sphere such that it and its subgraphs have complements which are homeomorphic to the connect-sum of handlebodies.    We relate the freeness index to the long-standing orientable cycle double cover conjecture (OCDCC).  In order to make this connection, we exploit the fact that once a graph is embedded into a closed orientable surface $F$, we can then embed $F$ into $S^3$.  If both the embedding of the graph into the surface, and of the surface into the 3-sphere, are sufficiently ``nice'', we are able to conclude that the induced embedding of the graph into the 3-sphere satisfies some strong freeness conditions.     We use this to show that a cubic graph satisfying the OCDCC has freeness index at least 2.      While we conjecture that every connected cubic graph has freeness index at least 2,  a connected cubic graph with freeness index 1 would provide a counter-example to the OCDCC.    

Here is an outline of the paper:

In Section \ref{tech}, we give essential definitions and a technical claim from Heegaard theory. 

In Section \ref{0-free}, as a warm-up, we show every graph can be embedded in the 3-sphere with handlebody complement, i.e., has freeness index at least 0. 

In Section \ref{panelled}, we recall work of Robertson, Seymour and Thomas and relate their notion of a ``panelled" embedding to the freeness index. 

In Section \ref{K_6} we explicitly calculate the freeness index of $K_6$ (showing it is exactly 8) and of the Petersen graph (showing that it is exactly 4). 

In Section \ref{index 1} we show that every graph has freeness index at least 1. 

In Section \ref{other} we show that graphs that have 2-factors with few components, including the flower snarks,  have freeness index greater than one.

In Section \ref{strong} we relate the freeness index of $\Gamma$ to the existence of a strong or polyhedral embedding of $\Gamma$ into a closed orientable surface.   Note that for cubic graphs the existence of such a strong embedding is equivalent to the graph having an orientable cycle double cover.

In Section \ref{dual}, we describe a dual formulation of the freeness index, which allows us to invoke the forbidden minors theorem of Robertson and Seymour.

In Section \ref{ocdcc}, we show that a cubic graph with a strong embedding into a closed orientable surface (equivalently, satisfying the orientable cycle double cover conjecture) has a 2-free embedding.   

In Section \ref{question}, we close with a question and conjecture.

\section{Technical details}
\label{tech}

Our graphs will always be finite. \\

\noindent Notation: If $X$ and $Y$ are topological spaces and $\phi$ is an embedding of $X$ into $Y$, we denote the image of $X$ in $Y$ by $X_\phi$.

We  denote $S^3$ with a small open neighborhood of $\Gamma_{\phi}$ removed by $S^3-\Gamma_{\phi}$.

\begin{Def} 
Let $\phi$ be an embedding of a finite graph $\Gamma$ into $S^3$, with $\Gamma_{\phi}$  the image of $\Gamma$ under $\phi$.   We say that $\Gamma_{\phi}$ is {\em free} (or {\em 0-free}) if $\pi_{1}(S^3-\Gamma_{\phi})$ is a free group.  This is equivalent to the  assumption that $(S^3-\Gamma_{\phi})$ is a connect-sum of handlebodies.  
\end{Def}

The embedding ${\phi}$ induces an embedding of any subgraph $\Gamma'$ of $\Gamma$.   We denote by $\Gamma'_{\phi}$ the image of a subgraph $\Gamma'$ of $\Gamma$ under the induced embedding.  

\begin{Def}
We say the embedding $\Gamma_{\phi}$  is {\em k-free} if $\Gamma'_{\phi}$ is free for all subgraphs $\Gamma'\subset{\Gamma}$ obtained from $\Gamma$ by removing the interiors of at most $k$ edges.     
\end{Def}

\begin{Def}
The {\em freeness index} of $\Gamma$ is the maximal $k$ for which $\Gamma$ admits a $k$-free embedding.    
\end{Def}

\begin{Def}
Given a graph $\Gamma$, a {\em cycle} in $\Gamma$ is a closed path in $\Gamma$. 

A {\em cycle double cover} of $\Gamma$ is a collection $\mathcal{C}=\{C_1, . . . , C_k\}$ of cycles in $\Gamma$ such that each edge of $\Gamma$ is visited by precisely two cycles $C_i,C_j$, $i\neq{j}$.
\end{Def}

\begin{Def}
Given a graph $\Gamma$, a {\em 2-factor} of $\Gamma$ is a collection of disjoint cycles in $\Gamma$ which contain every vertex.   As an example, a graph with a Hamiltonian cycle has a 2-factor (the Hamiltonian cycle). 
\end{Def}

\begin{Def}
A  {\em Heegaard surface for $S^3$} is a closed orientable surface $F$ embedded in $S^3$ such that $F$ splits $S^3$ into two handelbodies, $H_1$ and $H_2$.   There exists a Heegaard surface for $S^3$ of every genus. 
\end{Def}

We will use the following well-known fact from 3-manifold theory (see \cite{H}, \cite{J}):\\

\begin{Claim}\label{Heegaard}

Let $F$ be a Heegaard surface for $S^3$ splitting $S^3$ into two handlebodies, $H_1$ and $H_2$.     Suppose there exists a pair of disks $D_1, D_2$, such that:
\begin{enumerate}
\item $D_i$ is properly imbedded in $H_i$.
\item $\partial{D_i}$ is essential in $F$.
\item $\partial{D_1}$ intersects $\partial{D_2}$ transversely in a single point.   

\end{enumerate}

\noindent Then the surface $F'$ obtained by compressing $F$ along $D_1$ (or $D_2$)  is still a Heegaard surface for $S^3$.  In particular both complementary pieces of $S^3-F'$ are handlebodies.     

\end{Claim}

\begin{Def} 
Such a pair of disks $D_1, D_2$ is called a {\em destabilizing pair} for $F$.   Reversing a destabilization, that is, taking the connect sum of a Heegaard surface $F$ with an unknotted torus in $S^3$, is called a {\em stabilization} of $F$, and results in a new Heegaard surface of one greater genus.     
\end{Def}  

\begin{Coro}

Let $\Gamma_{\phi}$ be an embedding of a connected graph $\Gamma$ into $S^3$ such that $S^3-\Gamma_{\phi}$ is a handlebody $H_1$.    Suppose there exists a compressing disk $D_1$ of $H_1$ such that $\partial(D_1)$ runs over an edge $e$ of $\Gamma_{\phi}$ exactly once.     Then $S^3-(\Gamma_{\phi}-e)$ is a handlebody $H_{1}'$.   \\

\end{Coro}

\noindent Proof:
Let $H_{2}=N(\Gamma)$ be a neighborhood of $\Gamma$, and let $m\subset{\partial{H_2}=\partial{H_1}}$ be a meridian of $e$, bounding a disk $D_2$ in $H_{2}$.  The disks $D_1,D_2$ are a destabilizing pair for the Heegaard splitting of $S^3$ defined by $F=\partial{H_1}=\partial{H_2}$, and we apply Claim \ref{Heegaard}. \\

\section{All graphs have $0$-free embeddings}
\label{0-free}

It is easy to show that every connected graph $\Gamma$ admits an embedding $\theta$ into a closed surface $F$ such that:

 \noindent (1) the connected components of $F-\Gamma_{\theta}$ are all open discs, and \\
 (2) the surface $F$ is orientable,\
 
 by taking a small tubular neighborhood of $\Gamma$, projecting a copy of $\Gamma$ onto the bounding surface $F'$, and then surgering $F'$ along essential simple closed curves disjoint from the image of $\Gamma$ as much as possible to obtain $F$.

Given an embedding $\theta$ of $\Gamma$ into $F$ satisfying (1) and (2),  we may then further embed $F$ via a map $\phi$  into $S^3$ as a Heegaard surface, i.e., such that $F_\phi$ separates $S^3$ into two handlebodies.

 Given such an $F_\phi$, since the complementary components of $\Gamma_\theta$ in $F$ are discs, the complement of $\Gamma_{\phi(\theta)}$ in $S^3$ is also a handlebody, constructed by attaching together the inner and outer handlebodies defined by $F_\phi$ along the disk faces of $F_\phi-\Gamma_{\phi(\theta)}$.

 Hence $\pi_{1}(S^3-\Gamma_{\phi(\theta)})$ is a free group, and so this is a 0-free embedding of $\Gamma$ into $S^3$.

\section{Freeness index and panelled embeddings}
\label{panelled}

Some graphs with $n$ edges have freeness index $n$, i.e., they can be embedded in $S^3$ so that the fundamental group of the complement of every subgraph is free.    These graphs were studied by Robertson, Seymour and Thomas in \cite{RST}, from whom we adopt some concepts and notation.     

\begin{Def} An embedding $\phi$ of a graph $\Gamma$ into $S^3$ is called {\em panelled} if every cycle in the embedded graph bounds a disk with interior disjoint from the rest of the graph.  
\end{Def}

In \cite{RST} it is shown that an embedding $\phi$ is panelled if and only if  $\pi_{1}(S^3-\Gamma'_{\phi})$ is free for every subgraph $\Gamma' \subset{\Gamma}$. Every planar graph admits a panelled embedding, simply by embedding it first into a 2-sphere and then embedding the 2-sphere into $S^3$.  The class of graphs that admit panelled embeddings is larger than the planar graphs; the non-planar complete graph on $5$ vertices, $K_5$, for example, also admits a panelled embedding.    However not all graphs admit such an embedding; some graphs, such as $K_6$, are known to be {\em intrinsically linked}, i.e., for any embedding $\phi$ into the 3-sphere, there will always be a subgraph $\Gamma'$ such that $\pi_{1}(S^3-\Gamma'_{\phi})$ is not free \cite{CG}\cite{S}.  In \cite{RST} it is shown that  $\Gamma$ admits a panelled embedding if and only if none of the seven graphs of the Petersen family is contained in $\Gamma$ as a minor. 

The freeness index refines the class of graphs that do not admit panelled embeddings.

\section{Freeness index of $K_6$ and the Petersen graph}
\label{K_6}

From the previous discussion it follows that $K_6$ has freeness index strictly less than 15 (its number of edges), and that any graph that fails to admit a panelled embedding has freeness index strictly smaller than its number of edges.\\

We will show:\\
 \begin{Prop}
 \label{proposition}
 
 \begin{enumerate}
\item The freeness index of $K_6$ is 8.
\item The freeness index of the Petersen graph is 4. 
\end{enumerate}
\end {Prop} 

This means that there exists an embedding of $K_6$ into $S^3$ for which the complement of the graph is free, and such that when we remove (up to) any 8 edges from the graph, the complement remains free, but there is no such embedding for 9 edges.   Similarly, there exists an embedding of the Petersen graph in $S^3$ such that the complement is free and when we remove (up to) any 4 edges from the graph the complement remains free, but there is no such embedding for 5 edges.     

We will need the following technical lemma:

\begin{Lem}\label{Lemma1} 
Let $\Gamma^{\prime}$ be a connected planar graph, and let $\Gamma$ be obtained from $\Gamma'$ by attaching a single additional edge.   
Then $\Gamma$ has a panelled embedding.
\end{Lem}

\noindent Proof of Lemma \ref{Lemma1}:\\
\noindent Let $\phi'$ be a planar  embedding of $\Gamma'$ into $S^2\subset{S^3}$.   $S^2$ divides $S^3$ into two 3-balls; we will think of one as lying ``above'' $S^2$ and one as lying ``below''. 

 Let $P$ be a path in $\Gamma'_{\phi'}$ connecting the endpoints of $e$.    Define an embedding $\phi$ of $\Gamma$ by attaching the edge $e$ to a parallel push-off of $P$ disjoint from $S^2$.    We can assume $e$ lies slightly “above” the $S^2$ containing $\Gamma'_{\phi'}$.   Any cycle completely contained in $\Gamma'_{\phi'}$ bounds a disk slightly below  $S^2$ with interior disjoint from the rest of $\Gamma_{\phi}$. 

Assume $C$ is a cycle containing the edge $e$. Then $C$ consists of a path $Q$ in $S^2$ together with the edge $e$, and $C$ bounds a disk lying slightly above $S^2$ with interior disjoint from the rest of $\Gamma_{\phi}$. 

Hence every cycle of $\Gamma_{\phi}$ bounds a disk with interior disjoint from bounds a disk with interior disjoint from $\Gamma_{\phi}$.

\noindent Proof of Proposition \ref{proposition}:

$(1)$ Conway and Gordon \cite{CG} and Sachs \cite{S} showed that any embedding of $K_6$ into $R^3$ contains a non-trivial link, consisting of two cycles of length 3, obtained by removing 9 edges from $K_6$. Thus for any embedding of $K_6$, there is always some choice of 9 edges such that removing them gives a non-freely embedded subgraph.   Hence $K_6$ is not 9-free.

\begin{figure}[htbp] %  figure placement: here, top, bottom, or page
   \centering
   \includegraphics[width=3in]{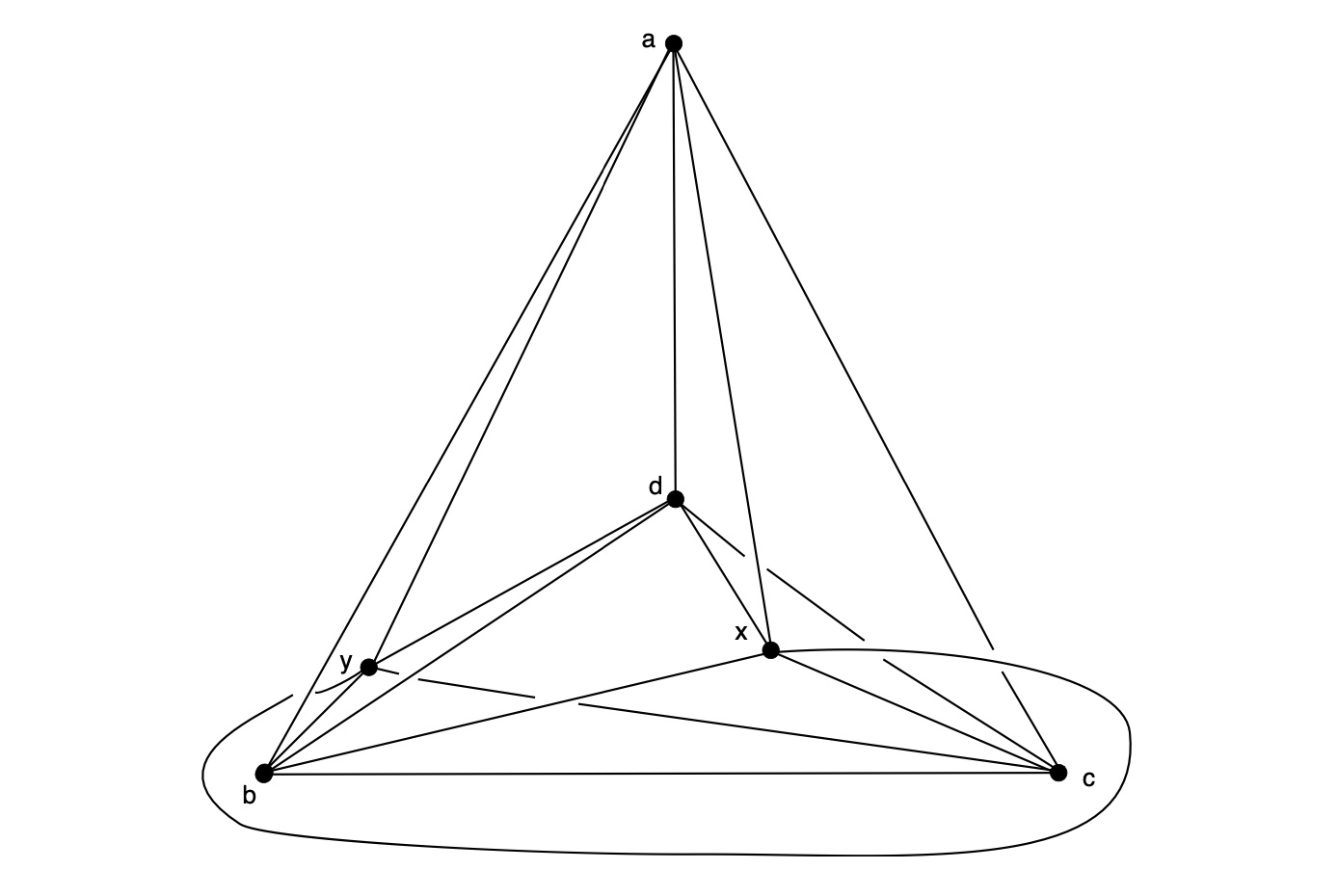} 
   \caption{an embedding of $K_6$}
   \label{fig1}
\end{figure}

\begin{figure}[htbp] %  figure placement: here, top, bottom, or page
   \centering
   \includegraphics[width=3in]{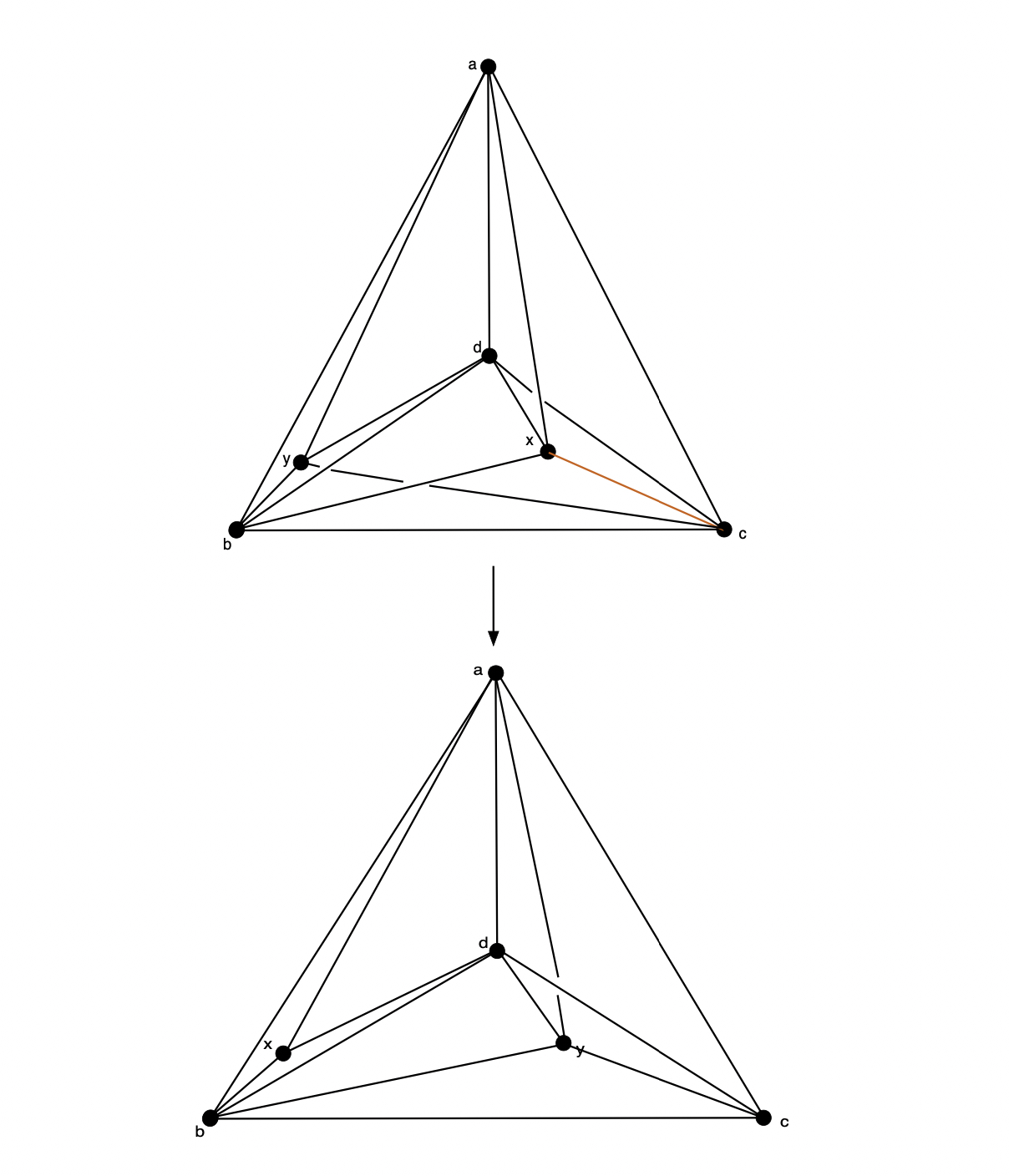} 
   \caption{removing $xc$ and isotoping}
   \label{fig2}
\end{figure}

\begin{figure}[htbp] %  figure placement: here, top, bottom, or page
   \centering
   \includegraphics[width=3in]{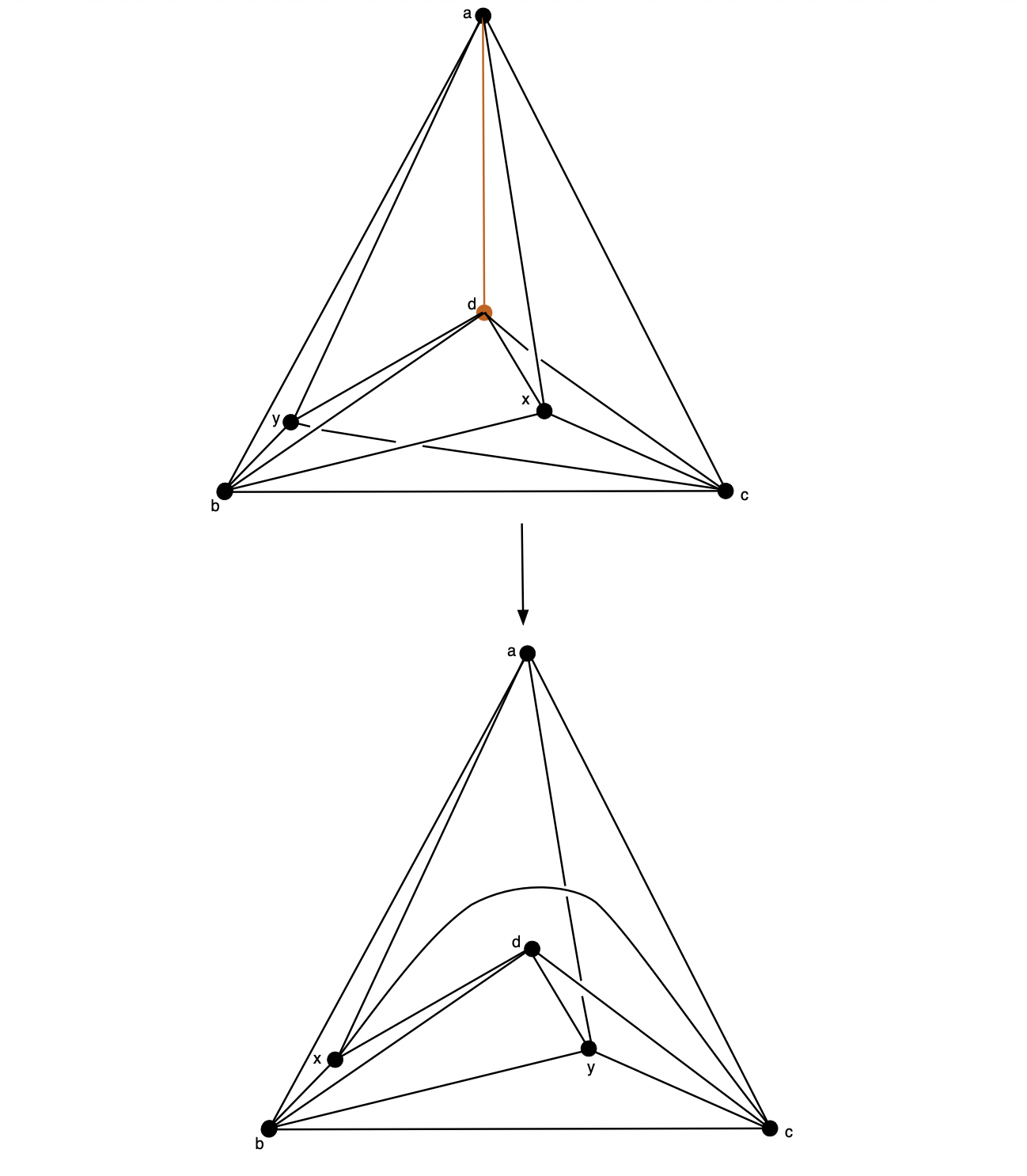} 
   \caption{removing $ad$ and isotoping}
   \label{fig3}
\end{figure}

We now show that the embedding of $K_6$ shown in Figure \ref{fig1} is 8-free.

In the pictured embedding, the unique non-trivial link is the Hopf link given by the pair of 3-cycles $C_1 = (xdy)$ and $C_2 = (abc)$. Observe that every other 3-cycle in this embedding bounds a disk disjoint from the rest of the graph, hence $C_1\cup{C_2}$ is the unique nontrivial sublink in this embedding. We now show that removing any subset of up to eight edges from this embedding results in a graph with free complement.

Suppose we remove up to eight of the edges that connect $C_1$ to $C_2$. The complement of the Hopf link formed by $C_1\cup{C_2}$ is homeomorphic to $T^2\times{I}$ , and the remaining connecting edges (of which there must be at least one) can be isotoped to be parallel and monotonic in the $I$ factor of this product. The complement of this is free.

Suppose next that at least one edge of $C_1$ or $C_2$ is among the edges removed. By symmetry we can assume this edge is $xy$. Call the resulting (embedded) graph  $\Gamma'$. If $xy$ is the only edge removed, the complement is free by observation (see the top half of Figure \ref{fig2}). In fact this embedding is panelled. It is possible to see this by removing any additional edge, then isotoping the resulting graph, and then applying Lemma \ref{Lemma1}. An additional removed edge is either one of $xc,xd,xa,xb,ya,yb,yc,yd$ or one of $ad,ab,ac,bc,bd,cd$.  In Figure \ref{fig2} we illustrate the isotopy in the case that $xc$ is removed, and in Figure \ref{fig3} we illustrate the isotopy in the case that $ad$ is removed. It is helpful to picture the vertex $x$ lying slightly above   a plane containing  $a,b,c,d$,  and $y$ lying slightly below, at the start of the isotopy. The other cases are similar.

\begin{figure}[htbp] %  figure placement: here, top, bottom, or page
   \centering
   \includegraphics[width=3in]{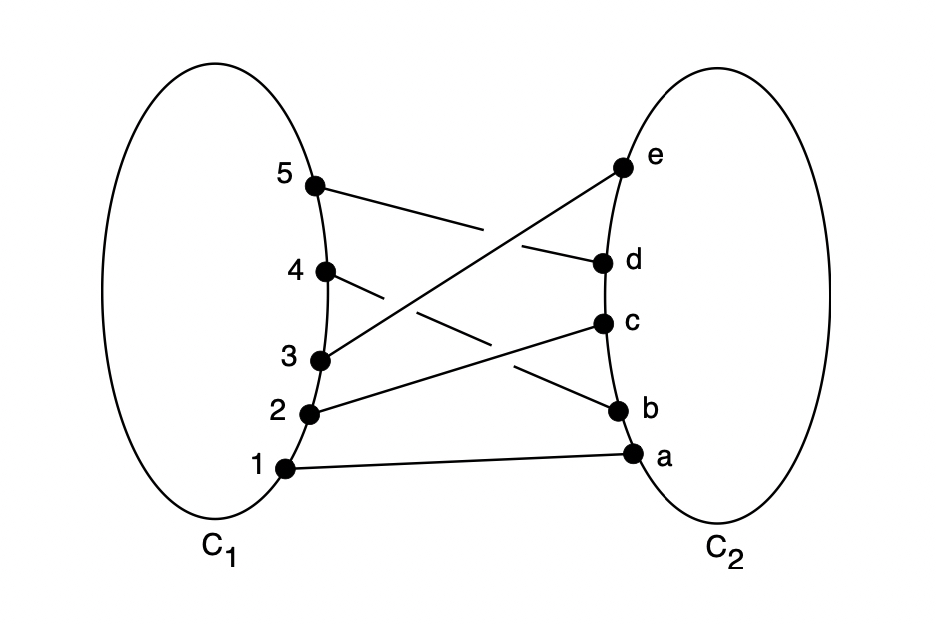} 
   \caption{an embedding of the Petersen graph}
   \label{fig4}
\end{figure}

$(2)$ Robertson, Seymour and Thomas \cite{RST} showed that any embedding of the Petersen graph into $R^3$ must contain a non-trivial link. Any such link must consist of two cycles of length 5 (the length of the shortest cycle in the Petersen graph).  Hence no embedding of the Petersen graph is 5-free. We show that the embedding of the Petersen graph shown in Figure \ref{fig4} is 4-free.

\begin{figure}[htbp] %  figure placement: here, top, bottom, or page
   \centering
   \includegraphics[width=3in]{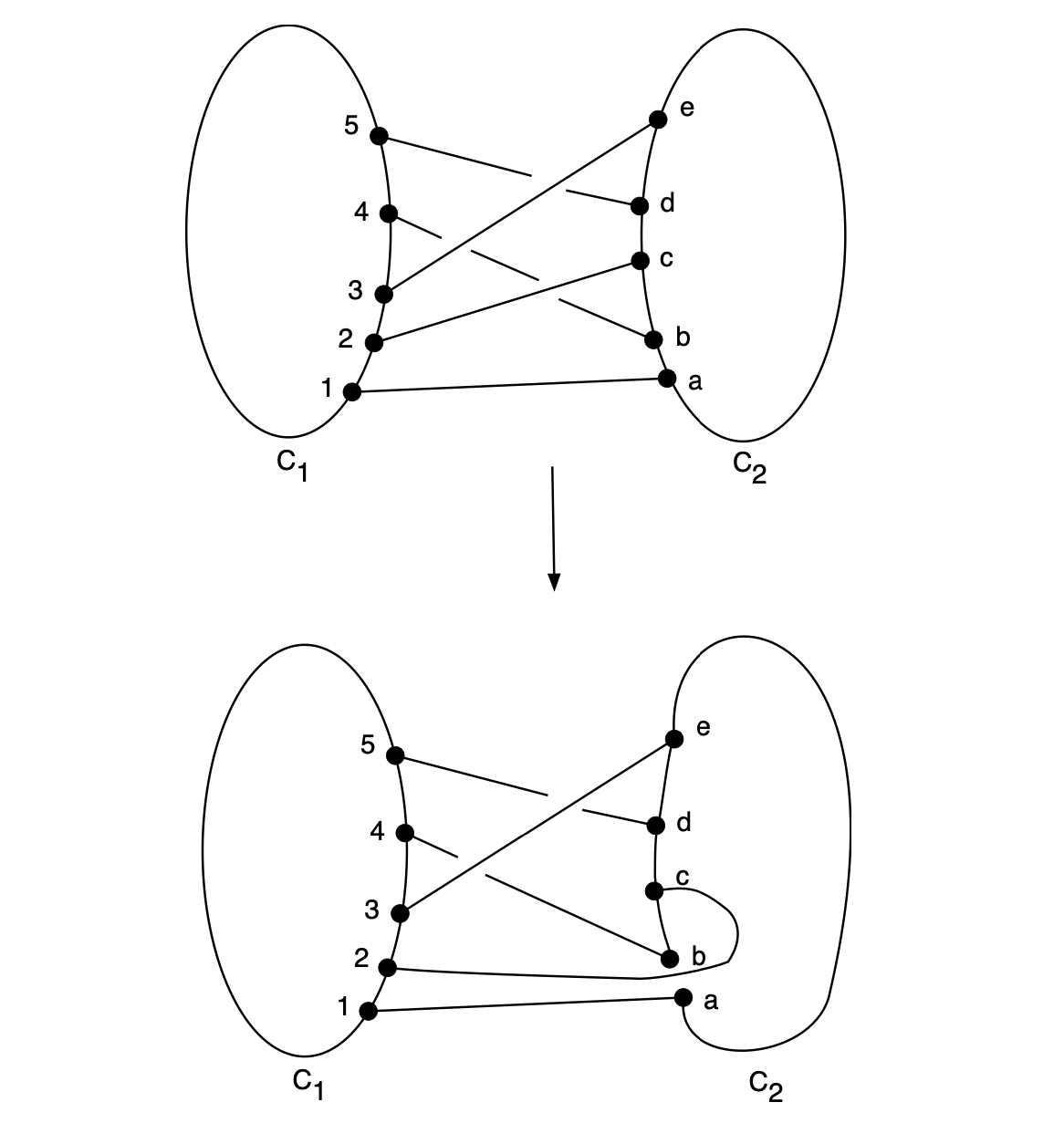} 
   \caption{removing the edge $ab$ from the Petersen graph}
   \label{fig5}
\end{figure}

This pictured embedding is not 5-free since it contains a Hopf link, consisting of the cycles $H_1 = (5dc21)$ and $H_2 = (34bac)$. As in the case of the embedding of $K_6$ above, this is the unique non-trivial link in this embedding. The argument is similar; any link must consist of two disjoint 5-cycles. We check directly that all 5-cycles except for $H_1$ and $H_2$ bound disks disjoint from the rest of the graph. This is obviously true for the cycles $C_1$ and $C_2$. Any other 5-cycle must consist of two edges connecting $C_1$ to $C_2$, and either two edges contained in $C_1$ and one in $C_2$ or vice versa. In fact, picking any two edges connecting $C_1$ to $C_2$ uniquely defines a 5-cycle. It is thus easy to check that each of these except $H_1$ and $H_2$ bounds a disk disjoint from the rest of the graph.

Suppose we remove up to four edges connecting $C_1$ to $C_2$. The complement of the two disks bounded by $C_1$ snd $C_2$ is isotopic to $S^2\times{I}$, and the remaining connecting edges (and there must be at least one) can be isotoped to be monotonic in the $I$ factor. The complement of any such graph is free.

Suppose we remove up to four edges connecting $H_1$ to $H_2$. As in the $K_6$ case, the complement of the Hopf link formed by $H_1$ and $H_2$ is isotopic to $T^2\times{I}$, and the remaining connecting edges (there must be at least one) can be isotoped to be monotonic and parallel in the $I$ factor. The complement of such a graph is again free.

Suppose at least one edge is removed that is contained in one of the $C_i$'s and contained in one of the $H_i$'s. Removing any such edge results in a subgraph satisfying the hypotheses of Lemma \ref{Lemma1}, hence is panelled. We illustrate this case in Figure \ref{fig5}.

We have shown that every graph has a  0-free embedding but not all graphs have a 5-free embedding. In the next section we show that we can improve the lower bound for the freeness index of any graph by proving that every graph has a 1-free embedding.

\section{All graphs have $1$-free embeddings}
\label{index 1}

In fact we can do better than 1-free; we show that every graph has a 1-free embedding with any selection of disjoint cycles embedding as the unlink:

\begin{Thm}{\label{1-free}}
Let $\Gamma$ be a connected graph, and let $\mathcal{C}$ be any collection of disjoint cycles in $\Gamma$. 
Then $\Gamma$ has a 1-free embedding $\phi$ in $S^3$. 
Moreover $\phi$ can be chosen so that $\mathcal{C}_\phi$
is the unlink and the components of $\mathcal{C}_\phi$ bound a collection of disjoint disks with interiors disjoint from $\Gamma_\phi$.
\end{Thm}

\noindent Proof:

Assume the statement is false and suppose  $\Gamma$ is a minimal counterexample to the theorem;  
that is, assume any graph with fewer edges than $\Gamma$ satisfies the theorem. 
Note $\Gamma$ is not a single cycle. Let $\mathcal{C}$ be a (possibly empty) 
collection of disjoint cycles in $\Gamma$. Since $\Gamma$ is connected and 
is not a single cycle, we can choose an edge $e$ in $\Gamma$ which
 is not contained in any cycle in $\mathcal{C}$. Let $u,v$ be the endpoints of $e$.  Let $\Gamma'=\Gamma-e$ denote $\Gamma$ with the interior of $e$ removed.  

Suppose $e$ is a bridge in $\Gamma$. Since $\Gamma$ is a minimal counterexample, the components of $\Gamma'$ each has a 1-free embedding in $S^3$ with the induced cycles from $\mathcal{C}$ satisfying the theorem. Use these 1-free embeddings to construct a 1-free embedding of their union by embedding each into a 3-ball in $S^3$ separated by a 2-sphere $S$.  Attach the embedded components of $\Gamma'$ by an edge $e$ intersecting the sphere $S$ in a single point and disjoint from the disks bounded by the cycles in $\mathcal{C}$.   

 This is a 1-free embedding of $\Gamma$;  this follows from a variant of the classic light-bulb trick, for details see \cite{HT}.

Assume $e$ is a not a bridge. Since $\Gamma$ is a minimal counterexample and $\Gamma'$ is connected, $\Gamma'$ has a 1-free embedding $\phi'$ in $S^3$ with the cycles from $\mathcal{C}$ satisfying the theorem.  We denote this embedding $\Gamma_{\phi}'$.
 Let $\mathcal{D}$ be the collection of disks bounded by the components of $\mathcal{C}$. Let $P$ be a path in $\Gamma'$ from $u$ to $v$. 

\begin{figure}[htbp] %  figure placement: here, top, bottom, or page
   \centering
   \includegraphics[width=3in]{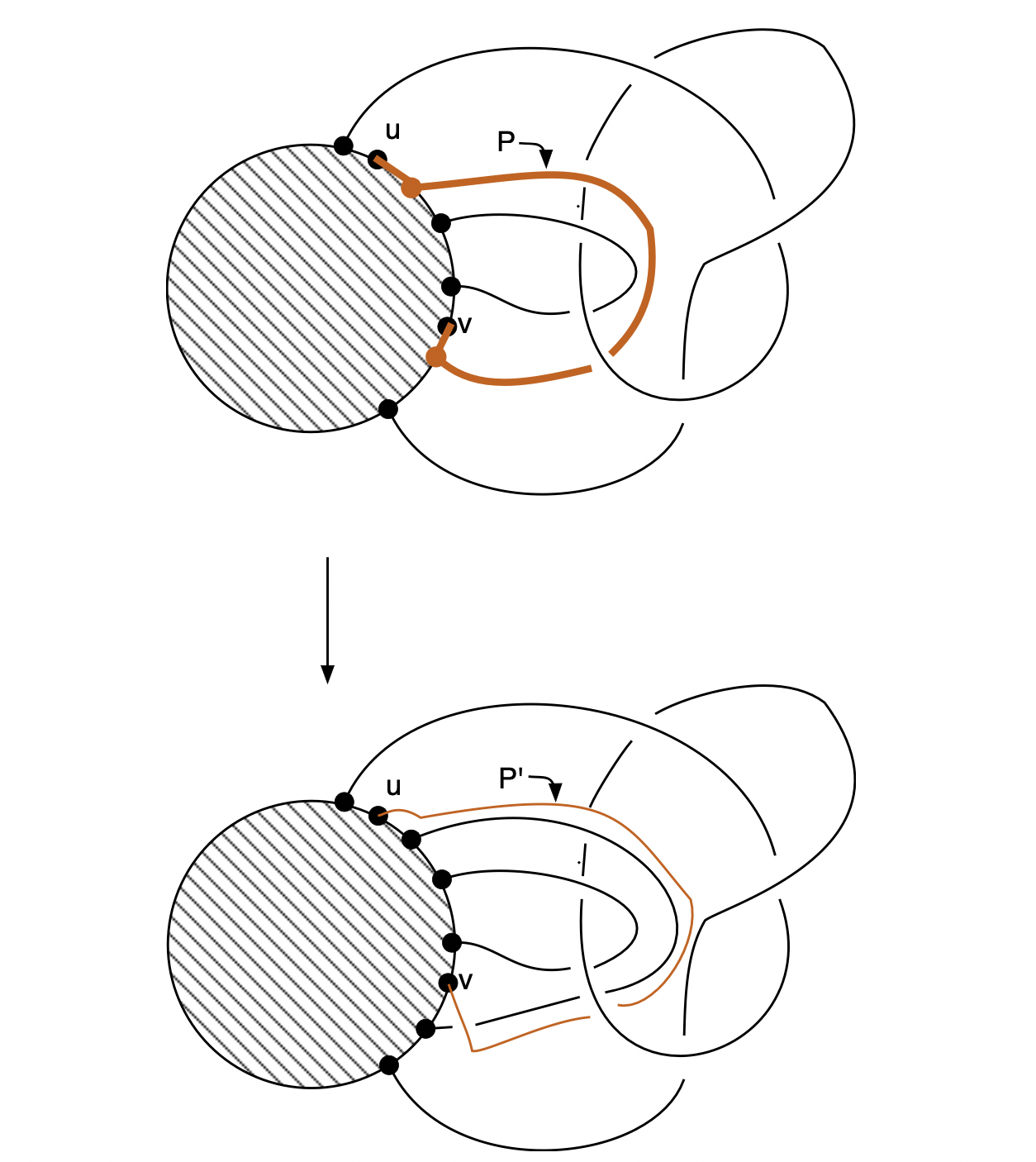} 
   \caption{attaching the edge e to $\Gamma'_\phi$}
   \label{fig6}
\end{figure}

\begin{Claim}
We can choose a push-off $P'$ of $P_\phi$ in $S^3$  that is disjoint from the interiors of $\mathcal{D}$.
\end{Claim}
Proof. An example of a 1-free embedding of a $\Gamma'$ is shown in Figure \ref{fig6}, with a single disk $D$ shown shaded. The path $P_\phi$ is shown in red in the top of the figure, with its push-off $P'$ shown in the bottom. 

\bigskip

Then $\Gamma_{\phi}'\cup{P'}$ is an embedding of $\Gamma$ into $S^3$, where $e$ is mapped to $P'$. We denote this embedding by $\Gamma_\phi$.   By construction, the cycles of $\mathcal{C}$ form an unlink in   $\Gamma_\phi$ and the components of $\mathcal{C}$ bound disjoint disks with interiors disjoint from  $\Gamma_\phi$. We now show $\Gamma_\phi$ is 1-free.

The cycle $E = P\cup{P'}$ bounds a disk $G$ with interior disjoint from the rest of  $\Gamma_\phi$. This suffices to show that $\Gamma_\phi$ is free. 

If we remove an edge from $\Gamma_\phi$  that is disjoint from $E$, the complement remains free; this follows from the 1-free property of $\Gamma_{\phi}'$ and the fact that $e$ is added to $\Gamma_{\phi}'$ along a boundary parallel arc, with the parallelism given by $G$. If we remove $e$ itself the complement is free since $\Gamma_{\phi}'$ is free. Now suppose we remove an edge $f$ of P. Since the disk  $G$ runs over $f$ once, by Claim \ref{Heegaard} it follows that $\Gamma_{\phi}-f$ is still free.

\section{Freeness index for some classes of graphs}
\label{other}

In this section we prove a more general version of  Theorem \ref{1-free}  for certain classes of graphs.

\begin{Thm}\label{Hamiltonian} 
If $\Gamma$ has a Hamiltonian cycle $C$, then $\Gamma$ has a 3-free imbedding ${\phi}$.  Moreover ${\phi}$ can be chosen so that $C_\phi$ bounds a disk $D$ with the interior of $D$ disjoint from the edges of $\Gamma_{\phi}$.

\end{Thm}

\begin{Coro}
\label{complete}
The complete graph on $n$ vertices has a 3-free imbedding.   
\end{Coro}

\noindent Proof of Corollary \ref{complete}:\\
The complete graph on $n$ vertices has a Hamiltonian cycle. 

\begin{Thm}\label{2-factor,2}
Assume $\Gamma$  is cubic and connected.    If $\Gamma$ has a 2-factor with two components,  $C_1, C_2$ then  it has a 2-free imbedding ${\phi}$. Moreover ${\phi}$ can be chosen so that ${C_1}_\phi, {C_2}_\phi$ bound a pair of disjoint disks $D_1, D_2$ with the interior of $D_i$ disjoint from the edges of $\Gamma_{\phi}$.

\end{Thm}

\begin{Thm}\label{2-factor,3}
Assume $\Gamma$  is cubic and 3-connected.   If $\Gamma$  has a 2-factor with three components, it has a 2-free imbedding.  Moreover the embedding can be chosen so that the image of the 2-factor is an unlink bounding disks with interiors disjoint from the rest of the graph. 
\end{Thm}

\begin{Thm}\label{2-factor,4}
Assume $\Gamma$  is cubic and 3-connected.   If $\Gamma$  has a 2-factor with four components, it has a 2-free imbedding.  Moreover the embedding can be chosen so that the image of the 2-factor is an unlink bounding disks with interiors disjoint from the image of all other edges of $\Gamma$. 

\end{Thm}

Proof of Theorem \ref{Hamiltonian}:\\

We proceed by induction on $n$ where $n$ is the number of chords of $C$.    

The statement is obviously true for $n=1,2$.

Assume the statement is true for $n-1$, $n>2$ chords.

Let $\Gamma$ be a graph with Hamiltonian cycle $C$ and $n$ chords.     Let $e$ be a chord of  $\Gamma$.     Let $\Gamma'$ be obtained from $\Gamma$  by deleting the interior of $e$.   $C$ is still a Hamiltonian cycle of $\Gamma'$.    By the inductive hypothesis, there is an imbedding $\phi'$ of  $\Gamma'$ which is 3-free and in which $C_{{\phi}'}$ bounds a disk $D$ with interior disjoint from  $\Gamma'_{\phi'}$.   Extend this to an imbedding $\phi$ of $\Gamma$ by imbedding the edge $e$ into the disk $D$; $e$ divides the disk $D$ into two subdisks $D_1$ and $D_2$.   

To avoid too many subscripts, we denote  specific edges of $\Gamma$ and their images in $S^3$ with the same letters.

 By pushing $e$ slightly off of $D$, we can see that $C_{{\phi}}$ still bounds a disk (which we continue to call $D$) with the interior of $D$ disjoint from the edges of $\Gamma_{\phi}$.    We need to show that this imbedding $\Gamma_{\phi}$ is 3-free.     Let $f,g,h$ be any three edges of $\Gamma$.    If one of these three is actually the edge $e$, removing the (images of the) three edges from $\Gamma_{\phi}$ gives the same complement as removing two edges from  $\Gamma'_{\phi'}$.   Since $\Gamma'_{\phi'}$ is 3-free, the complement is free.    
 
Now assume none of the edges $f,g,h$ is $e$.   We can assume that at most one of these edges, say $f$, is contained in $\partial{D_1}$.   \\
 
 Claim: The complement of $(\Gamma_{\phi})-(g\cup{h})$  is free.\\
 Proof of Claim: The complement of  $\Gamma'_{\phi'}-(g\cup{h})$ is free since $\Gamma'_{\phi'}$ is 3-free.    Since $e$ is parallel to a path in $\Gamma'_{\phi'}$ via the disk $D_1$, the complement of $(\Gamma_{\phi})-(g\cup{h})$ is also free.  \\
 
 Suppose $f$ is not an edge of $D_1$.  Since  $\Gamma'_{\phi'}-(f\cup{g}\cup{h})$ is free, and $e$ is parallel to a path in $\Gamma'_{\phi'}-(f\cup{g}\cup{h})$ via $D_1$, it follows that $\Gamma_{\phi}-(f\cup{g}\cup{h})$ is free.    
 
 If $f$ is an edge of $D_1$, Claim \ref{Heegaard} implies that $\Gamma_{\phi}-(f\cup{g}\cup{h})$ is free.\\

\noindent Proof of Theorem \ref{2-factor,2}:\\

We induct on the number $n$ of edges in  ($\Gamma$ -(edges of $C_1\cup{C_2}$)).    

If any edge of $\Gamma$ is a chord in $C_1$, say, we can remove that edge and use the inductive hypothesis together with the arguments from the proof of Theorem \ref{Hamiltonian} to obtain the desired conclusion. 

\begin{figure}[htbp] %  figure placement: here, top, bottom, or page
   \centering
   \includegraphics[width=3in]{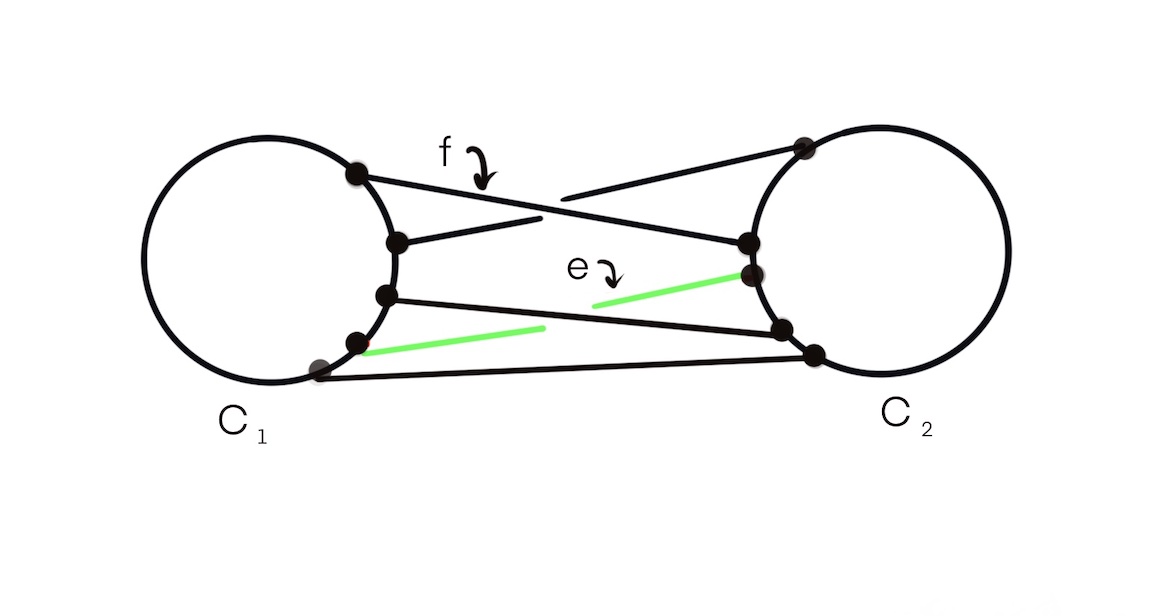} 
   \caption{two component 2-factor}
   \label{fig7}
\end{figure}

\begin{figure}[htbp] %  figure placement: here, top, bottom, or page
   \centering
   \includegraphics[width=3in]{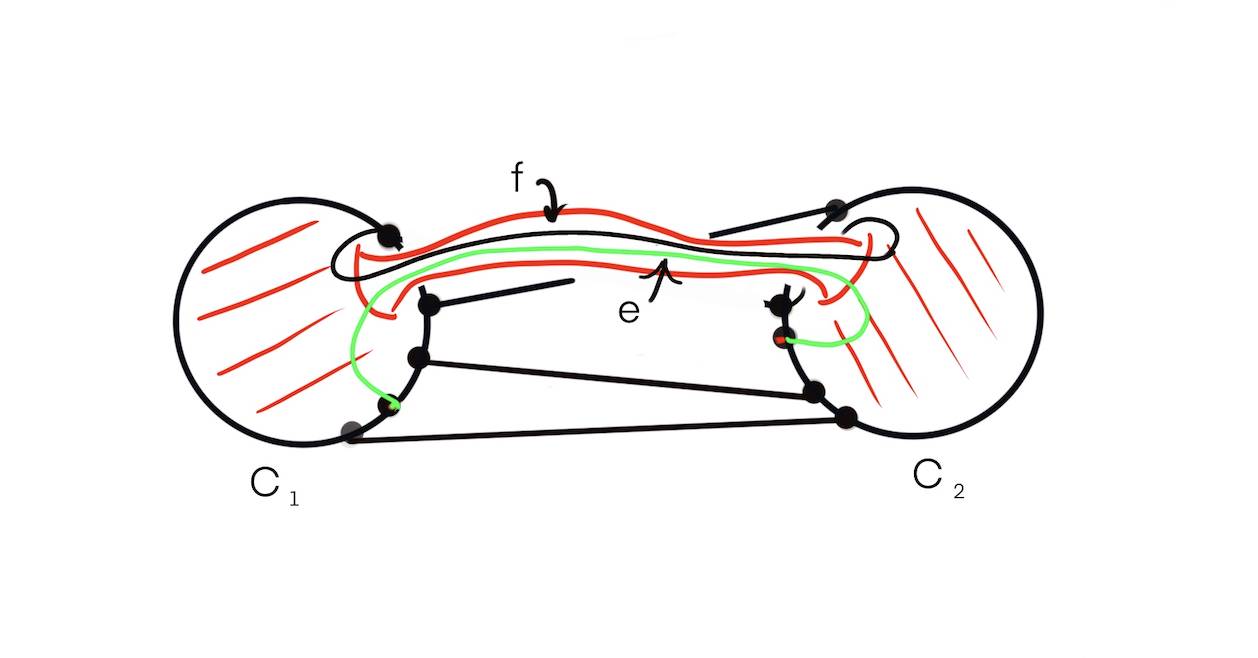} 
   \caption{}
   \label{fig8}
\end{figure}

 If no edge of $\Gamma$ is a chord in $C_1$ or in $C_2$, then there are at least two edges $e,f$ connecting $C_1$ to $C_2$ (else the Theorem is trivially true).   Remove $e$ from $\Gamma$ and amalgamate its endpoints to obtain $\Gamma'$. Apply the inductive hypothesis to obtain an imbedding $\phi'$ of $\Gamma'$, where the cycles $C_1,C_2$ of the induced 2-factor bound disjoint disks  $D_1, D_2$ which are also disjoint from the rest of the edges of $ \Gamma'_{\phi'}$.   Let $E$ be a disk constructed by banding $D_1$ and $D_2$ together along a band following $f$ (see Figures \ref{fig7} and \ref{fig8}).  The boundary of $E$ runs twice over the edge $f$.  Note this can be done so that $E$ is disjoint from $D_1\cup{D_2}$ and also from the other edges of $\Gamma'_{\phi'}$. 
Now imbed $e$ into $E$, to obtain an imbedding $\phi$ of $\Gamma$ (Figure \ref{fig8}).   $e$ divides  $E$ into two subdisks $E_1$ and $E_2$.    Since $e$ is parallel to a path (indeed to two paths) in  $ \Gamma'_{\phi'}$, it follows that $\Gamma_{\phi}$ is free.    We need to show  $\Gamma_{\phi}$ is 2-free.   

Let $g,h$ be any two edges of $\Gamma$. If one of $g,h$ is the edge $e$, the result follows by our hypothesis that $ \Gamma'_{\phi'}$ is 2-free.  We now assume neither $g$ nor $h$ is the edge $e$.    \\

Suppose in addition that neither $g$ nor $h$ is the edge $f$.    Then one of the subdisks $E_1$, $E_2$ has at most one edge removed from its boundary, and the result follows from the proposition together with the fact that $ \Gamma'_{\phi'}$ is 2-free. \\

It remains to consider the case where we remove the edge $f$ and we remove another edge $h\neq{e}$ from the imbedding. \\

\begin{Claim}  
$S^3-(\Gamma_{\phi} -(f\cup{h}))$ is homeomorphic to $S^3-(\Gamma'_{\phi'}-(h))$, and hence is free.   
\end{Claim}

\noindent Proof of Claim:    Either the boundary of $E_1$ or the boundary of $E_2$ (or perhaps both) does not contain $h$.    Suppose $\partial(E_1)$ does not contain $h$.    Using $E_1$ we can isotop $e$ (possibly sliding the endpoints of $e$ over vertices) to be parallel to $f$.  So $S^3-(\Gamma_{\phi}-(f\cup{h}))$ is homeomorphic to $S^3-(\Gamma'_{\phi'}-(h))$, as required.\\

\noindent  Proof of Theorems \ref{2-factor,3} and \ref{2-factor,4}:\\

A simple cubic graph with a 3- or 4-component 2-factor is either planar, or it has at least one chord with respect to a component of the 2-factor, or it has (at least) two edges connecting the same pair of components of the 2-factor.     Hence the claim can either be proved directly (in the planar cases) or by induction using the arguments in the proof of Theorem \ref{2-factor,2}.\\

\section{Freeness index and strong/polyhedral embeddings}
\label{strong}

For bridgeless cubic graphs $\Gamma$, while it is easy to find an embedding $\theta$ into a closed orientable surface $F$  such that every complementary region is a disk (Section \ref{0-free}), it may be difficult to find an embedding that is {\em strong}, i.e., so that the {\em closure} of every component of $F-\Gamma_{\theta}$ is an embedded disk.    Indeed proving that such an embedding always exists for this class of graphs  is equivalent to the orientable cycle double cover conjecture (OCDCC), an open problem.    An even stronger condition on graph embedding into a surface is as follows:
 
\begin{Def} Given a bridgeless cubic graph $\Gamma$  an embedding $\theta$ of $\Gamma$ into a closed orientable surface $F$ is  {\em polyhedral} if every complementary region of $F-\Gamma_{\theta}$ is a disk and any two such regions are either disjoint or share exactly one edge.    Equivalently, ${\theta}$ is polyhedral if every essential simple closed curve in $F$ intersects $\Gamma_{\theta}$ in at least three points.     
\end{Def}

If an embedding $\theta$ of $\Gamma$  in $F$  is strong, then for every edge $e$ of $\Gamma$, the complementary regions of $(\Gamma-e)_\theta$ in $F$ are still all discs.   As in Section \ref{0-free}, we can choose an embedding $\phi$ of $F$ so that $F$ is a Heegaard surface in $S^3$.    Then not only is $\pi_{1}(S^3-(\Gamma)_{\phi(\theta)})$ free,   $\pi_{1}(S^3-(\Gamma-e)_{\phi(\theta)})$ is also free for any edge $e$.  If the embedding of $\Gamma$  in $F$  is polyhedral, then a stronger result holds:

\begin{Thm}\label{polyhedral}
Assume $\Gamma$  is cubic and 3-connected. Let $\mathcal{C}=C_1,...C_k$ be any collection of disjoint cycles in $\Gamma$.   If $\Gamma$  has a polyhedral embedding into a closed orientable surface $F$, then it has a 2-free imbedding  into $S^3$.  Moreover the embedding can be chosen so that the image of $\mathcal{C}$ is the unlink bounding disks $D_1,...D_k$ with interiors disjoint from the images of all other edges of $\Gamma$.    

\end{Thm}

By considering the class of flower snarks, we can apply this to show:
  
\begin{Thm}\label{weaker}
The property of having a 2-free imbedding is strictly weaker than the property of having a polyhedral imbedding into a closed orientable surface. \\
\end{Thm}

Proof of Theorem \ref{polyhedral}:\\

Let ${\theta}$ be a polyhedral embedding of $\Gamma$ into the closed orientable surface $F$. $\Gamma_{\theta}$ divides $F$ into a collection of disks $D_1,...D_m$.    Imbed $F$ via $\phi$ as a Heegaard surface into $S^3$, splitting $S^3$ into the two handlebodies $H_1$ and $H_2$.    We can visualise the complement of the induced imbedding $\phi(\theta)$ of $\Gamma$ in $S^3$, denoted $\Gamma_{\phi(\theta)}$, as being constructed by attaching $H_1$ to $H_2$ along copies of the $D_{i}\times{I}$, i.e., $S^3-\Gamma_{\phi(\theta)}$ is constructed by gluing two handlebodies together along  1-handles, so  $S^3-\Gamma_{\phi(\theta)}$   is also a handlebody.      

Notice that for this construction to yield a handlebody complement for our imbedded graph, we only require that the complementary regions of the graph in the surface are disks.    Since  the imbedding of $\Gamma$ into $F$ is polyhedral, removing any two edges of $\Gamma$ from the imbedding still leaves all the complementary regions in $F$ as disks, hence the complement of the embedded graph with any two edges removed is still a handlebody.   

Now let $\mathcal{C}=C_1,...C_k$ be any collection of disjoint cycles in $\Gamma$, with polyhedral embedding $\theta$ into $F$.   We still need to find a 2-free embedding of $\Gamma$ into $S^3$ such that the image of $\mathcal{C}$ is the unlink bounding disks $D_1,...D_k$ with interiors disjoint from all other edges of the embedded $\Gamma$.    

$\mathcal{C}_{\theta}$ is some collection of simple closed curves on $F$ (some of which may be trivial, and some of which may be parallel).    As long as they are disjoint, however, we can find an embedding $\phi$ of $F$ as a Heegaard surface into $S^3$ such that each $(C_i)_{\phi(\theta)}$ bounds a disk in, say, $H_1$, and so the interior is disjoint from the edges of $\Gamma{\phi(\theta)}$ (see Figure \ref{fig11}).\\

\begin{figure}[htbp] %  figure placement: here, top, bottom, or page
   \centering
   \includegraphics[width=3in]{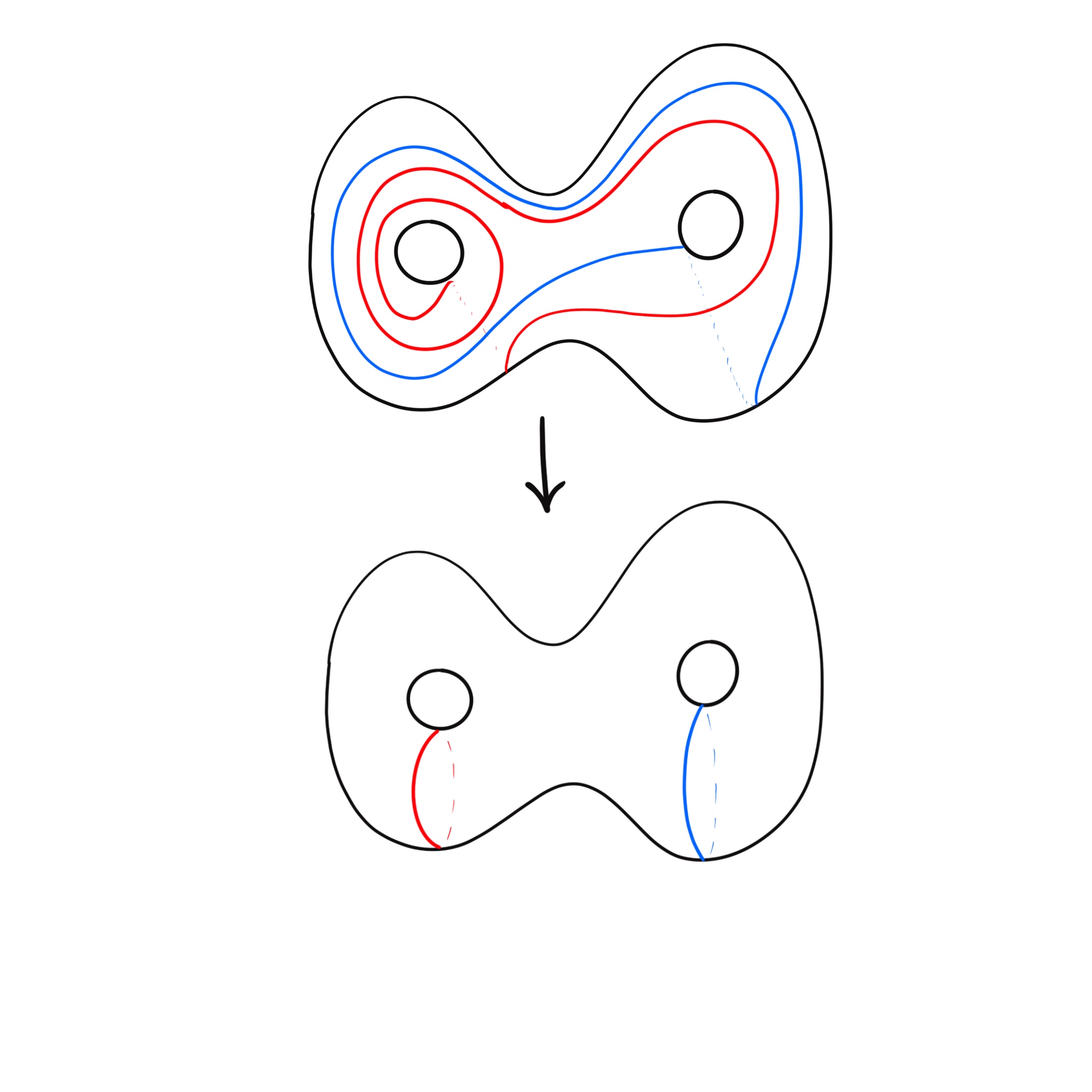} 
   \caption{re-embedding of a Heegaard surface with specified curves bounding disks in $H_1$}
   \label{fig11}
\end{figure}

\noindent Proof of Theorem \ref{weaker}:\\
The flower snarks have 2-component 2-factors, hence they have 2-free imbeddings by  Theorem \ref{2-factor,2}.  However they do not have polyhedral imbeddings into any closed orientable surface \cite{Sz}; indeed they don't have polyhedral imbeddings into any closed surfaces, orientable or not \cite{MV}.   \\

\begin{figure}[htbp] %  figure placement: here, top, bottom, or page
   \centering
   \includegraphics[width=3in]{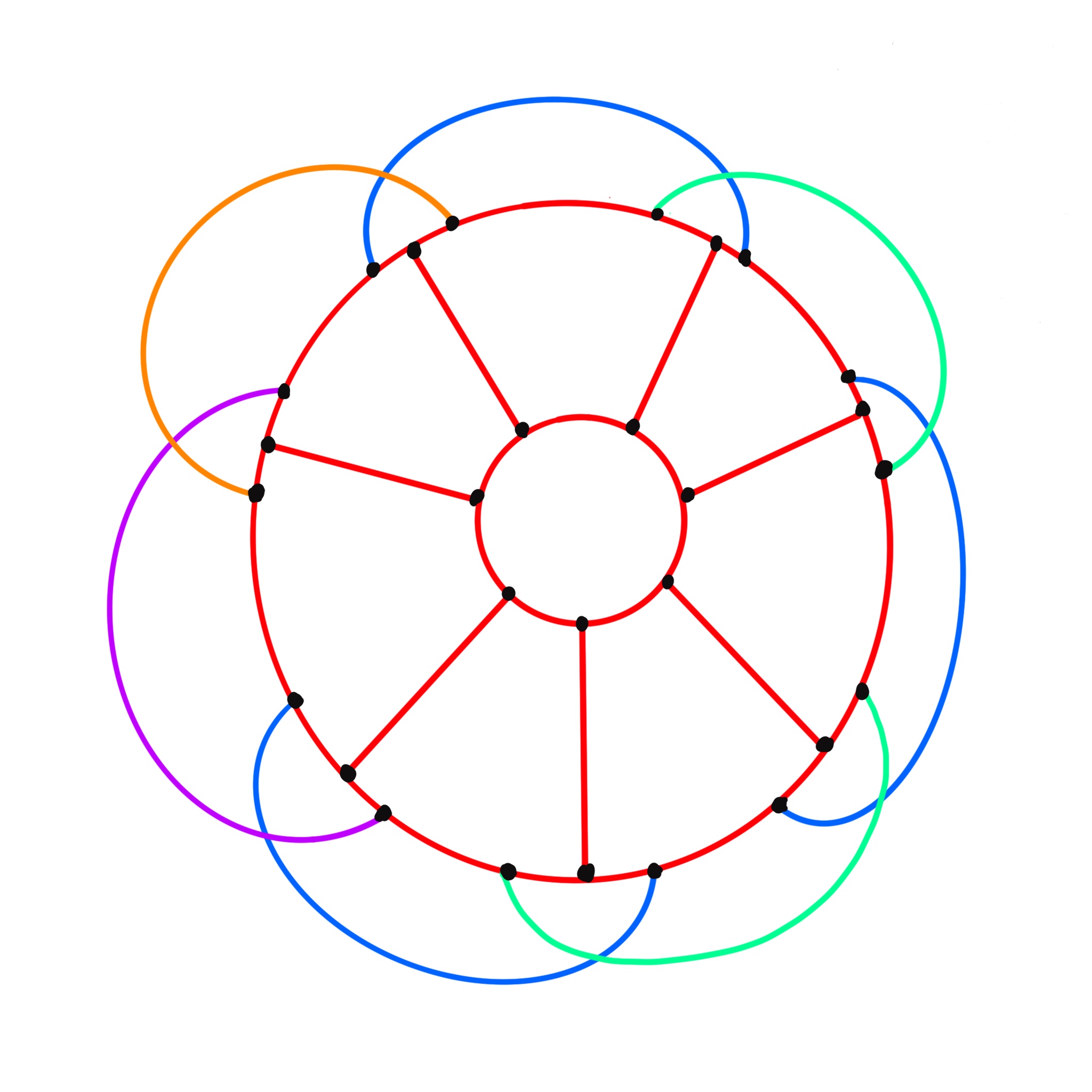} 
   \caption{$J_7$}
   \label{fig9}
\end{figure}

\begin{figure}[htbp] %  figure placement: here, top, bottom, or page
   \centering
   \includegraphics[width=3in]{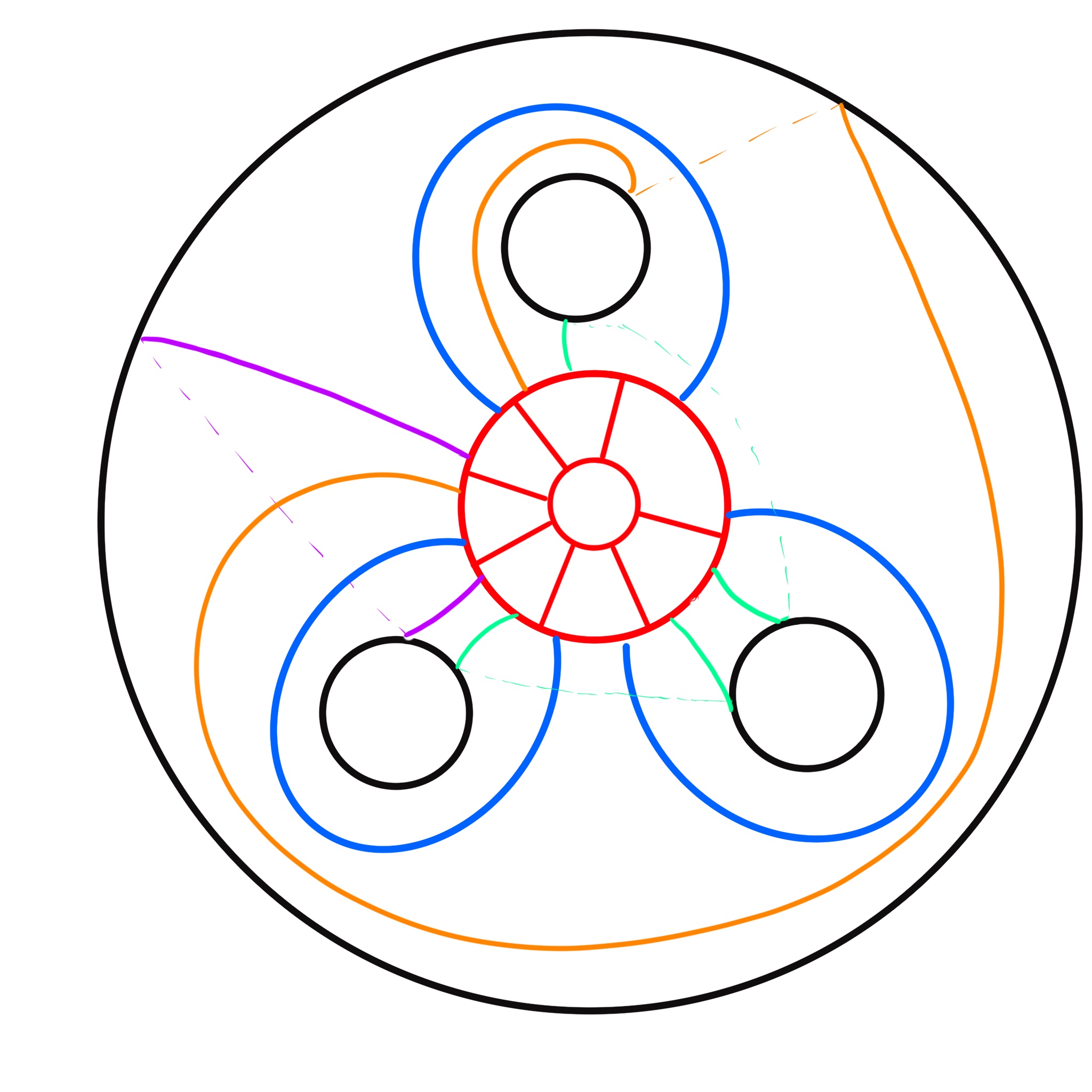} 
   \caption{a strong embedding of $J_7$ into a genus 3 surface}
   \label{fig10}
\end{figure}

While it is well-known that any potential counter-example to the Cycle Double Cover Conjecture must be a snark, the flower snarks do not provide such a counter-example, since they have strong embeddings into closed orientable surfaces.    We illustrate a strong embedding of the flower snark $J_7$ into a closed orientable genus 3 surface (see Figures \ref{fig9}, \ref{fig10} and \ref{fig11}), which generalizes easily to a strong embedding of $J_{2n+1}$ into a closed orientable surface $F$ of genus $n$ for any $n$.

\begin{figure}[htbp] %  figure placement: here, top, bottom, or page
   \centering
   \includegraphics[width=3in]{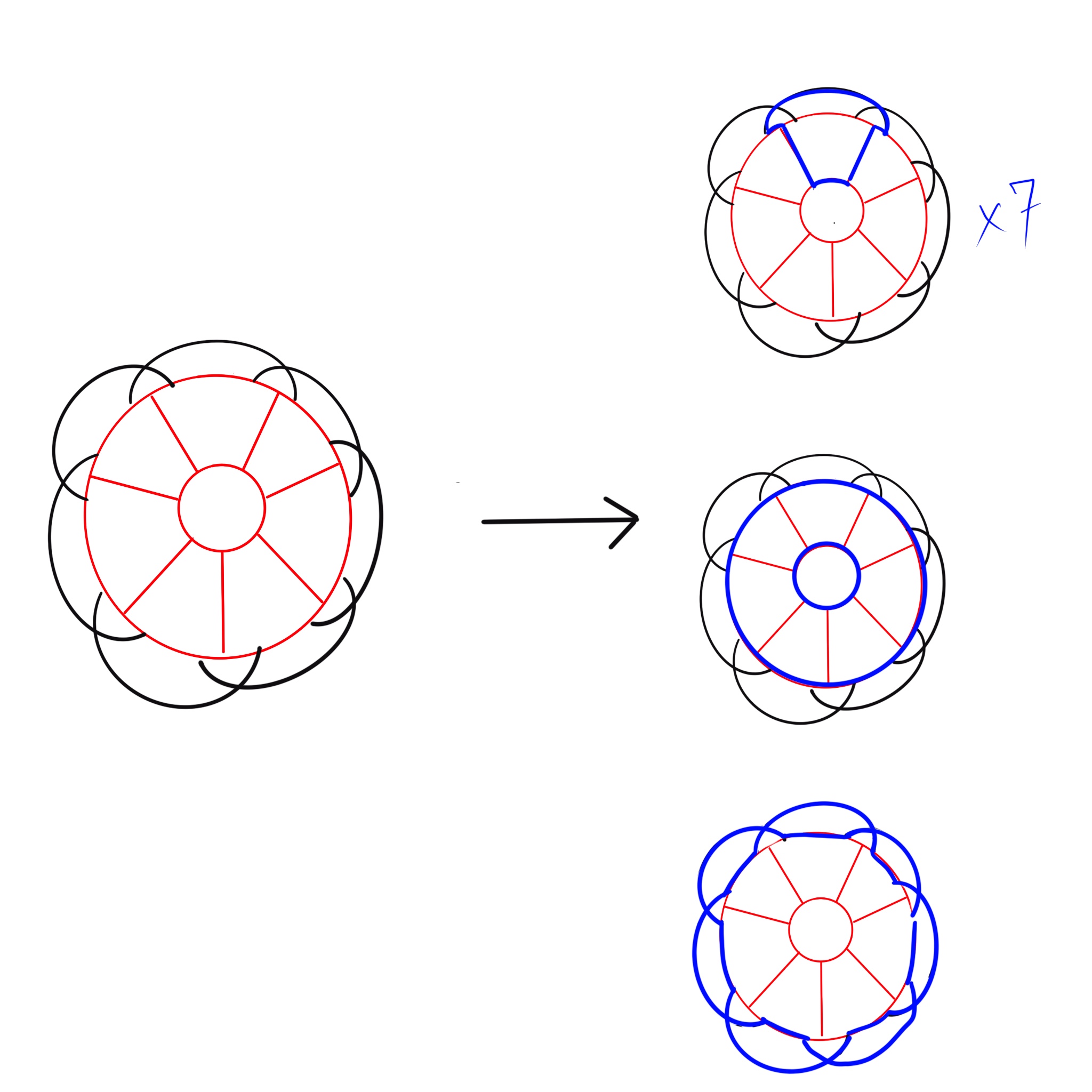} 
   \caption{strong embedding of $J_7$ into a genus 3 surface, described via an orientable cycle double cover}
   \label{fig11}
\end{figure}

\section{a dual property and forbidden minors}
\label{dual}

Thanks to Marc Lackenby for helpful comments on this section.

 We define a dual version of the freeness index:

\begin{Def} Let $\Gamma$ be a connected graph with $n$ edges, and $\phi$ an embedding of $\Gamma$ into the 3-sphere.    We say that $\phi$ is a {\em class $m$ embedding} if $S^3-{\Gamma}'$ is free for every subgraph $\Gamma'$ of $\Gamma$ with at least $m$ edges. 
\end{Def}

Note that if $\phi$ is a class $m$ embedding of $\Gamma$ it is also a class $p$ embedding for any $p$ with  $m\leq{p}$.   

\begin{Def}
The {\em class of $\Gamma$} is the smallest $m$ for which $\Gamma$ has a class $m$ embedding.   
\end{Def}

$\Gamma$  has a panelled embedding if and only if it has class $0$.

As additional examples, we have demonstrated that $K_6$ has class $7$, and the Petersen graph has class $11$.    

A re-statement of Theorem \ref{1-free} is that a graph $\Gamma$ with $n$ edges has class less than or equal to $n-1$.

We note that that the property of having a class $m$ embedding is closed under minors.   Hence, as for panelled embeddings, it follows immediately from \cite{RS} that there exists a finite collection of forbidden minors for each $m$:

\begin{Thm}\label{forbidden}
For each $m$ there is a finite collection of graphs $\Gamma_1,\Gamma_2,....,\Gamma_k$ such that if $\Gamma$ contains none of $\Gamma_1,\Gamma_2,....,\Gamma_k$ as a minor, then $\Gamma$ has a class $m$ embedding.  
\end{Thm}

\section{2-free embeddings and the orientable cycle double cover conjecture.}
\label{ocdcc}

The main result of this paper is the following: 
\begin{Thm}\label{main}
If $\Gamma$ has a strong embedding into a closed orientable surface $F$, then $\Gamma$ has a 2-free embedding into $S^3$.    
\end{Thm}

Proof:
The spirit of the proof is the same as for (the first part of) Theorem \ref{polyhedral}, but the details are trickier.     The difficulty with extending the argument about 2-freeness from a polyhedral to a strong embedding lies in the fact that given a strongly embedded graph in $F$, there may be many curves which intersect the graph essentially in just two points (see for example Figure \ref{fig10}).    When we remove  two edges corresponding to such a pair of points of intersection,  then when we construct the complement of the resulting graph in the 3-sphere, the inner- and outer- handlebodies defined by $F$  end up being attached together along an annulus and a collection of disks rather than simply along a collection of disks.     While it is possible to attach two handlebodies together along an annulus in such a way that the result is still a handlebody,  it must be done carefully.      This is the motivation behind  Lemma \ref{annuluslemma}. 

\begin{Lem}\label{annuluslemma}

Let $H_1$ and $H_2$ be two handlebodies, with $c_i$ a simple closed curve in $\partial{H_i}$.    Let $M$ be the quotient manifold  $H_1\cup_{c=c_1=c_2}H_2$.    Then $M$ is a handlebody or once-punctured handlebody if and only if either $c$ is primitive in $H_1$ or $H_2$, or $c$ is inessential in  $\partial{H_1}$ and in $\partial{H_2}$. 
\end{Lem}

This is a standard fact from Heegaard theory, with the slight modification that we consider the possibility that $c_i$ may be inessential in $\partial{H_i}$.    

The argument in one direction is straightforward; we show that if $M$ is a punctured handlebody then the conclusion holds.  \\

Proof of Lemma \ref{annuluslemma}:\\

Case 1:  $c_1$ does not bound a disk in ${H_1}$ and $c_2$ does not bound a disk in ${H_2}$.\\

Let $A$ be the properly embedded annular neighborhood of $c=c_1=c_2$ in $M$. Since $c_1$ does not bound a disk in ${H_1}$ and $c_2$ does not bound a disk in ${H_2}$, $A$ is incompressible in $M$.    Since the only incompressible and boundary incompressible surface properly embedded in M is a compressing disk, it follows that $A$ must be boundary compressible.  A boundary compressing disk for $A$ in $M$ arises exactly from a compressing disk for $H_1$ or $H_2$ intersecting $c_1$ or $c_2$ respectively in a single point.     \\

Case 2:  $c_1$ bounds a disk $D_1$ in ${H_1}$ and $c_2$ bounds a disk $D_2$ in ${H_2}$.\\
Then $D_1\cup{D_2}$ is an embedded 2-sphere in $M$, which must be boundary parallel.     Hence $c_1$ is inessential in $\partial{H_1}$ and $c_2$ is inessential in $\partial{H_2}$.\\

Case 3: $c_1$ bounds a disk $D_1$ in ${H_1}$ and $c_2$ does not bound a disk in ${H_2}$.\\
$M$ can be constructed by first attaching a 2-handle to $H_2$ along $c_2$ to obtain $M'$, and then attaching further 1-handles (corresponding to the 1-handles of $H_1$) to complete $M$, thus $M$ is a handlebody if and only if $M'$ is a handlebody.  But $M'$ is a handlebody if and only if $c_2$ is primitive in $H_2$ \cite{G}.

\begin{Thm}\label{2-points}
Let $\phi$ be a strong embedding of the graph $\Gamma$  into a closed orientable surface $F$.   Let $\mathcal{C}$ be the collection of all essential simple closed curves in $F$ which intersect $\Gamma_{\phi}$ in precisely two points.   Then $F$ can be embedded into $S^3$ so that: \\
1. $F$ is a Heegaard surface for $S^3$, splitting $S^3$ into two handlebodies $H_1$ and $H_2$,  and \\
2. for all $c\in{\mathcal{C}}$, either $c$ is primitive in $H_1$ or $H_2$,  or $c$ bounds a disk in  $H_1$ and $H_2$.
\end{Thm}

Proof:  

Let $\phi$ be a strong embedding of the graph $\Gamma$  into a closed orientable surface $F$.  We induct on the genus of $F$.   

If $F$ is a 2-sphere, the theorem is obviously true.   Note in this case that if $\Gamma$ is 3-connected, then the set $\mathcal{C}$ is empty.\\
Assume the theorem holds for genus$(F)=n-1$. \\

Assume genus$(F)=n\geq{1}$. \\

Case 1:   For every edge $e$ of $\Gamma$, there exists a $c\in{\mathcal{C}}$ such that $c$ intersects $e_{\phi}$. \\

\begin{Claim}\label{curves}
 In this case, there exists a pair of curves $c,d\in\mathcal{C}$ such that $c\cap{d}$ is one point, and that for all $f\in\mathcal{C}$, $f\cap{c}$ is empty or one point and $f\cap{d}$ is empty or 1 point.
 \end{Claim}
 
 We assume Claim \ref{curves} is true and complete the proof of Case 1. \\
 
 \begin{figure}[htbp] %  figure placement: here, top, bottom, or page
   \centering
   \includegraphics[width=3in]{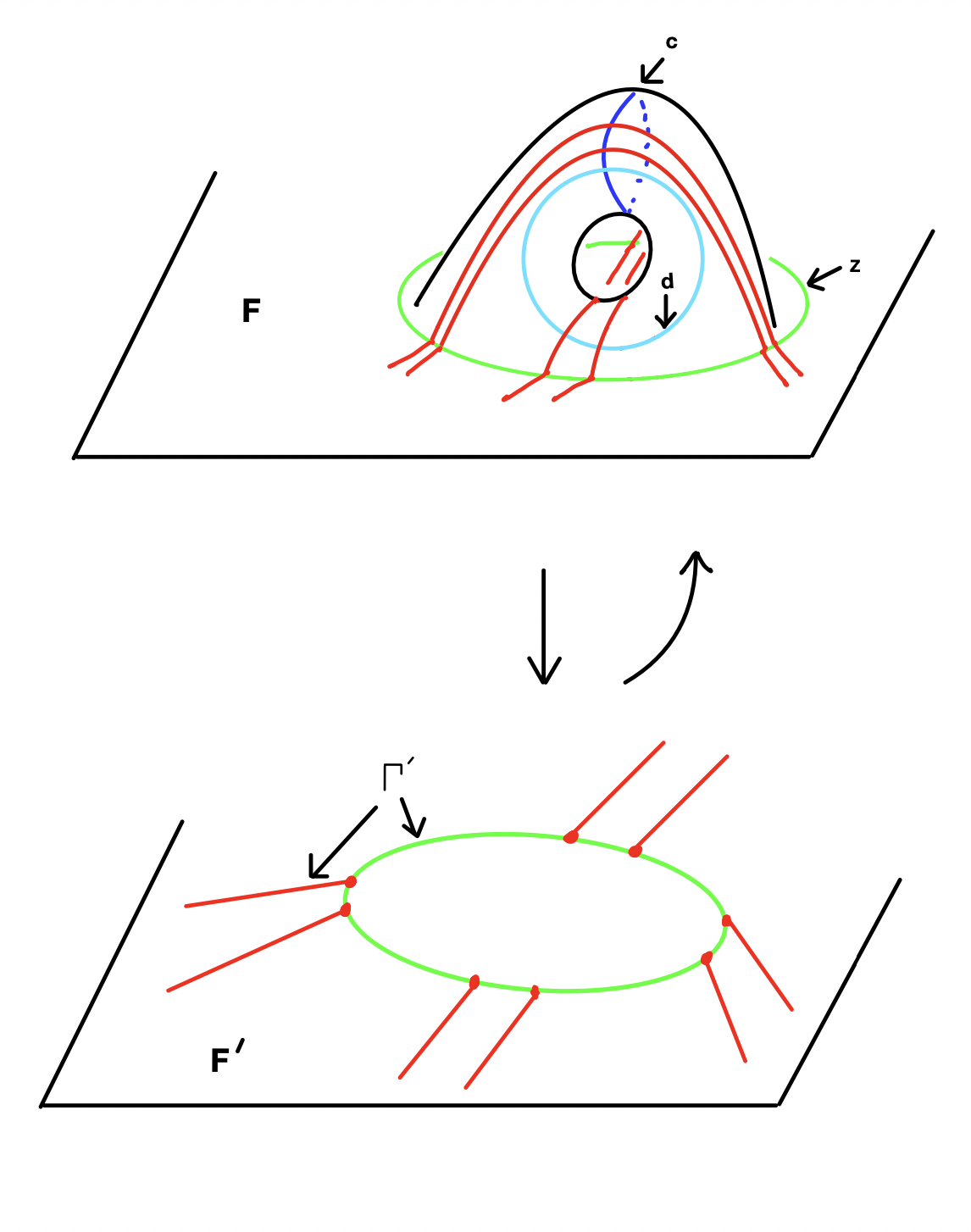} 
   \caption{two curves $c,d$ in $\mathcal{C}$ defining a punctured torus $T\in{F}$}
   \label{fig12}
\end{figure}

 Let $T$ be the punctured torus defined by the curves $c,d$, with $\partial{T}=z$.    We do the following modification of the pair $(F,\Gamma_\phi)$:\\
 Cut $F$ along $z$.    Replace $T$ with a disk to obtain a new surface $F'$  of one lower genus.    Replace $\Gamma_\phi$ by $\Gamma'=(\Gamma_\phi)-(\Gamma_{\phi}\cap{T})\cup({z})$ (see Figure \ref{fig12}, top to bottom).   Any simple closed curve in $F'$ which intersects $\Gamma'$ in one point must be disjoint from the disk bounded by $z$, which would would contradict that assumption that $\Gamma$ was strongly embedded in $F$.     Hence  $\Gamma'$ is strongly embedded in $F'$.     Since the genus of $F'$ is smaller than the genus of $F$, we can apply the inductive hypothesis to  the pair $(F',\Gamma')$ to obtain an embedding $\theta$ of $F'$ into $S^3$.     \\
 
 Now construct an embedding of $(F,\Gamma_\phi)$ by re-attaching $T$ to $F'_\theta$ so that the curves $c,d$ bound compressing disks on opposite sides of $F$ (see Figure \ref{fig12}, bottom to top).   This is a stabilization of $F'$, so this embedding of $F$ is still a Heegaard surface for $S^3$, splitting $S^3$ into two handlebodies $H_1$ and $H_2$.    Note that $c$ bounds a disk in, say, $H_1$ and is primitive in $H_2$,  and $d$ bounds a disk in $H_2$ and is primitive in $H_1$.   Any other  curve $f\in{\mathcal{C}}$  is either disjoint from $T\subset{F}$, hence satisfies Theorem \ref{2-points} by the inductive hypothesis, or it intersects (at least) one of $c,d$ in a single point and so is primitive, as required.

 Proof of Claim \ref{curves}:\\
 Because $\Gamma$ is strongly embedded in $F$, the components of $F-\Gamma_\phi$ are all disks $D_1,D_2,..,D_m$ with embedded closures in $F$.

\begin{figure}[htbp] %  figure placement: here, top, bottom, or page
   \centering
   \includegraphics[width=3in]{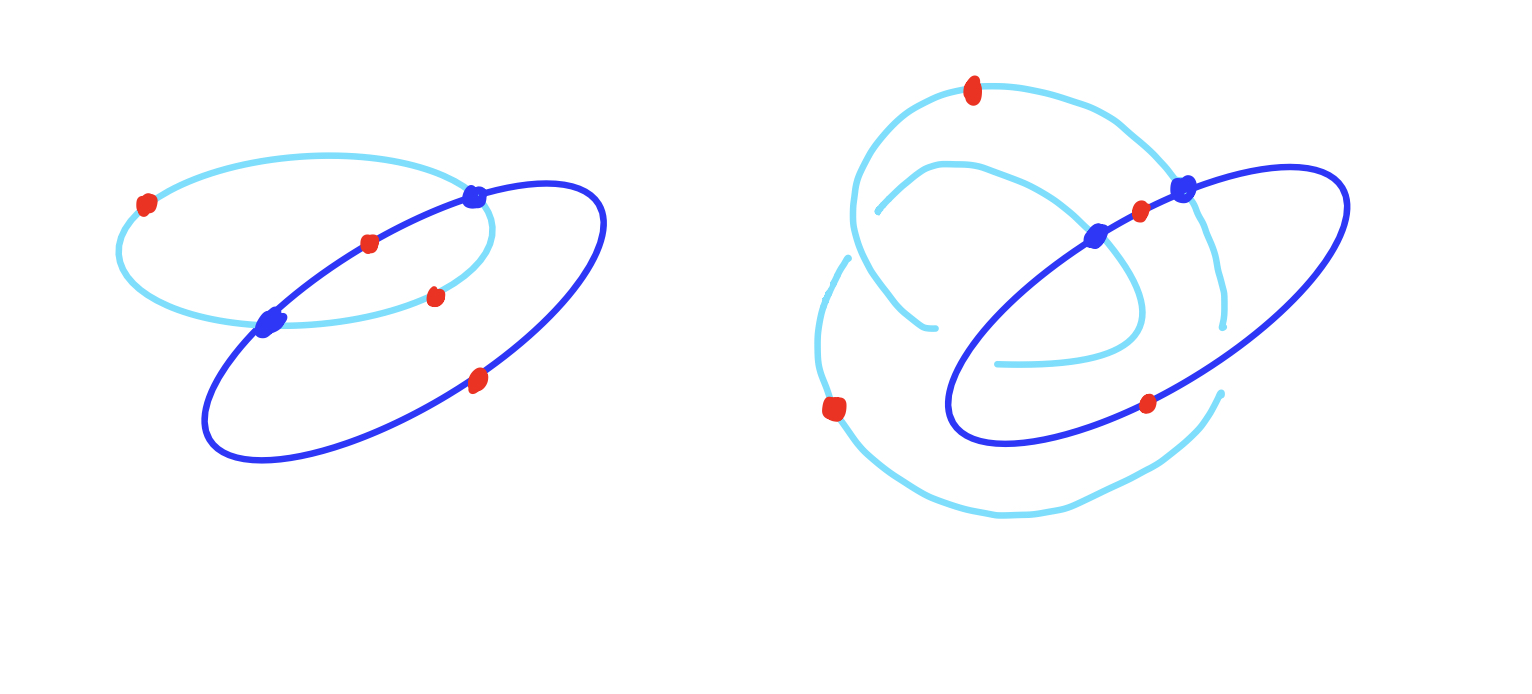} 
   \caption{two curves in $\mathcal{C}$; points of intersection with $\Gamma$ in red}
   \label{fig13}
\end{figure}

 We will need a few observations:\\
 1. for all $c,d\in\mathcal{C}$, $c$ intersects $d$ in 0,1 or 2 points.\\
 2. If $c,d\in\mathcal{C}$ intersect in 2 points, their points of intersection with $\Gamma_\phi$ must link the points of intersection of the two curves, as shown in Figure \ref{fig13}.\\
 3.  If $c,d\in\mathcal{C}$ intersect in 2 points $x,y$, we can (arbitrarily) ``smooth'' $x$ and then choose a smoothing of $y$ to obtain two new curves of $\mathcal{C}$. \\
 4. $c_1$ cannot connect adjacent edges in $D_1$.

 Proof of observation 1:    Suppose $c,d\in\mathcal{C}$, $c$ intersects $d$ in more than two points.     Then some arc $\alpha\subset{c}$ of $(c\cup{d}-(c\cap{d}))$ contains no points of intersection with $\Gamma_\phi$.    The endpoints of $\alpha$ divide $d$ into two arcs $d_1,d_2$, one of which, say $d_1$, intersects $\Gamma$ in 0 or 1 point.   Hence $\alpha\cup{d_1}$ is a simple closed curve on $F$ intersecting $\Gamma$ in 0 or 1 point.     If $\alpha\cup{d_1}$ is an essential curve, it contradicts the strong embedding of $\Gamma$.  If it is inessential, we could reduce the total number of intersection points in $\mathcal{C}$ by an isotopy.\\
 
 Proof of observation 2:   If not, then some arc $\alpha\subset{c}$ of $(c\cup{d}-(c\cap{d}))$ contains no points of intersection with $\Gamma_\phi$, and the argument proceeds as in the proof of observation 1.\\
 
\begin{figure}[htbp] %  figure placement: here, top, bottom, or page
   \centering
   \includegraphics[width=3in]{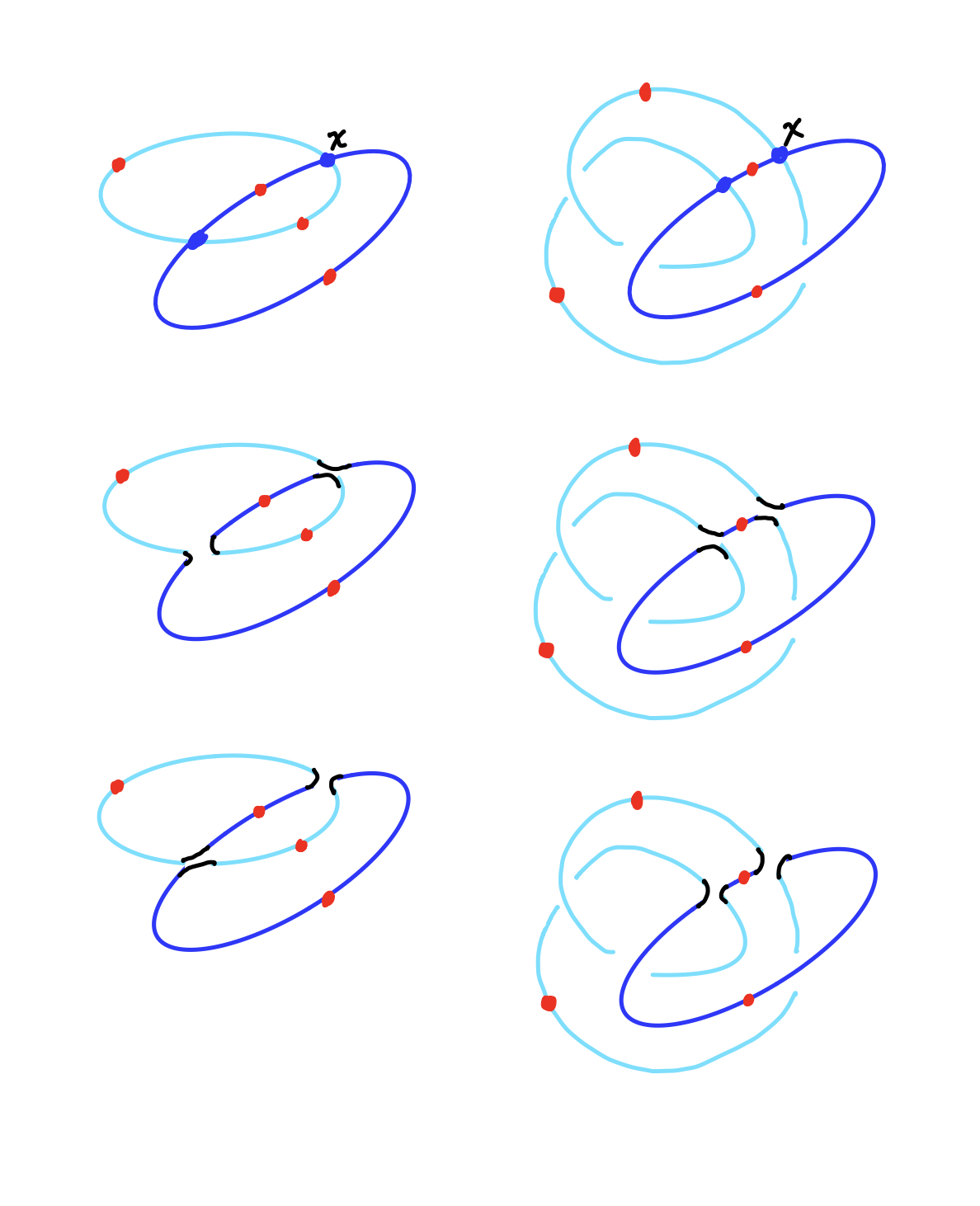} 
   \caption{arbitrary smoothing of a point $x$ of intersection}
   \label{fig14}
\end{figure}

 Proof of observation 3:   We use observation 2, and refer to Figure \ref{fig14}.\\
 
 We can describe $D_1$ as a polygon, with edges corresponding to edges of $\Gamma$.    \\
 Consider the collection of arcs of $D_1\cap{\mathcal{C}}$.     \\
 Let $c\in\mathcal{C}$.  $c$ intersects $D_1$ either in a single arc or not at all; suppose $c$ intersects $D_1$ in a single arc $c_1$.    We will say  $c_1$  {\em completes} to the curve $c$ in $F$.   If some $d\in{\mathcal{C}}$ intersects $c$ in two points, using observation 2 above, we can see that one of those points of intersection must lie in $D_1$.   \\
 
Proof of observation 4:    $c_1$ cannot connect adjacent edges in $D_1$.      If it does so, we can isotop $c$ across the vertex incident to both edges, and obtain a simple closed curve intersecting $\Gamma_\phi$ in a single point, a contradiction.

\begin{figure}[htbp] %  figure placement: here, top, bottom, or page
   \centering
   \includegraphics[width=3in]{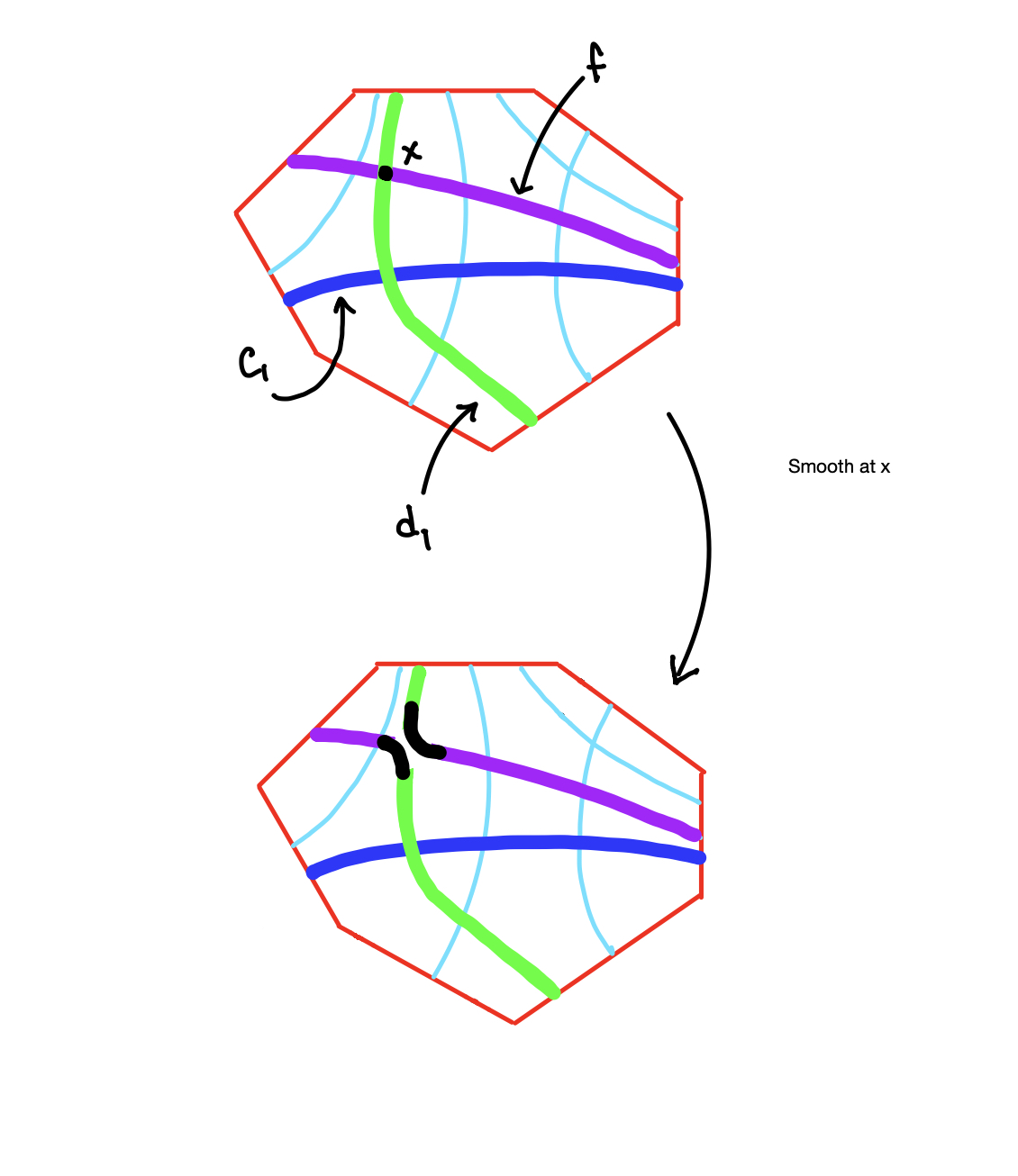} 
   \caption{Dark blue outermost arc $c_1$, green intersecting arc $d_1$ and purple $f\in\mathcal{C}$}
   \label{fig15}
\end{figure}

\begin{Def}  
$c_1$ is {\em outermost}  in $D_1$ if it divides $D_1$ into two disks, at least one of which does not completely contain any other arc of $\mathcal{C}\cap{D_1}$ (see Figure \ref{fig15} for an example of an outermost arc $c_1$.)
\end{Def}    

Subclaim:  If $c_1$ is outermost  in $D_1$, every curve $d\in{\mathcal{C}}$ intersects $c$ in 0 or 1 point, and there exists at least one $d\in{\mathcal{C}}$  which intersects in 1 point. 

Proof:    Assume $c_1$ is  outermost  in $D_1$.   Since $c_1$ cannot connect adjacent edges of $\partial{D_1}$, and since by hypothesis every edge of $\Gamma_\phi$ intersects some curve in $\mathcal{C}$, it follows that some curve $d\in{\mathcal{C}}$ must intersect $c_1$  in $D_1$.      Now suppose $d$ intersects $c$ twice.   Then we can smooth the points of intersection between $c$ and $d$ to obtain two new curves of $\mathcal{C}$.  But one of the arcs thus created contradicts the property that   $c_1$ was chosen to be outermost  in $D_1$.     \\

Subclaim:  We can find $d\in\mathcal{C}$ such that $d$ intersects $c$ (in one point) and for every $f\in{\mathcal{C}}$, $f$ intersects $d$ in 0 or 1 point.    \\

Proof:  Let $\mathcal{D}\in\mathcal{C}$ be the collection of all curves in $\mathcal{C}$ which intersect $c_1$ in $D_1$.    Let $d_1$ be outermost amongst all subarcs of curves of $\mathcal{D}$ in $D_1$, and let $d$ be the completion of $d_1$.      Suppose some $f\in\mathcal{C}$ intersects $d$ in two points $x,y$.    Then one of those points of intersection,  say $x$. must lie in $D_1$.    We can choose a smoothing of $x$ to produce another arc $\alpha$ which still intersects $c_1$ but which is contained in the outermost disk defined by $d_1$.     Given this smoothing at $x$, we can then smooth $y$ to produce two curves of $\mathcal{C}$, one of which contradicts our choice of $d_1$.     So no $f$  intersects $d$ in more than one point, as required.    This completes the proof of Case 1. \\

 Case 2:   Let $\phi$ be a strong embedding of the graph $\Gamma$  into a closed orientable surface $F$.  We now fix the genus of $F$, and induct on the number $n$ of edges of $\Gamma$ which do not intersect any curve in $\mathcal{C}$.\\
 
 $n=0$ is Case 1, above.     We assume the theorem holds for $n-1$ such edges.\\
 
 Suppose $\Gamma_\phi$ has $n\ge{1}$ edges that are disjoint from every curve of $\mathcal{C}$.   Let $e_\phi$ be one such edge.    Then $\phi$ induces a strong embedding of $\Gamma'=(\Gamma-e)$ into $F$ (by erasing the image of $e$ from $\Gamma_\phi$).   Apply the inductive hypothesis to $F$ to obtain an embedding $\theta$ of $F$ into $S^3$ so that all curves $\mathcal{D}$ intersecting $\Gamma'_\phi$ in two points are either primitive in $H_1$ or $H_2$ or bound a disk in both.  Now re-attach $e$.     Every $c\in{\mathcal{C}}$  was already  in $\mathcal{D}$, hence $\theta$  satisfies the conditions of the theorem for $\Gamma$.

  \bigskip

Proof of Theorem \ref{main} from Theorem \ref{2-points}:

Let $\phi$ be a strong embedding of the graph $\Gamma$  into a closed orientable surface $F$ and let $\mathcal{C}$ be the collection of all essential simple closed curves in $F$ which intersect $\Gamma_{\phi}$ in precisely two points.   By Theorem \ref{2-points} there exists an embedding $\theta$ of $F$ into $S^3$ such that $F$ is a Heegaard surface for $S^3$, splitting $S^3$ into two handlebodies $H_1$ and $H_2$,  and for all $c\in{\mathcal{C}}$, either $c$ is primitive in $H_1$ or $H_2$,  or $c$ bounds a disk in  $H_1$ and in $H_2$.   We will show that $\Gamma_{\theta(\phi)}$ is a 2-free embedding of $\Gamma$ into $S^3$.

Using the arguments from Theorem \ref{polyhedral}, it is easy to see that $\Gamma_\theta(\phi)$ is 1-free.    We need to consider what happens when we remove two edges from the embedded $\Gamma$.

If $e,f$ are a pair of edges of $\Gamma$ such that the regions of $F-(\Gamma-(e\cup{f}))_{\theta(\phi)}$ are all disks, then the arguments from Theorem \ref{polyhedral} show $(\Gamma-(e\cup{f}))_{\theta(\phi)}$ is still free.   Hence we only need to consider the case when  some $c\in{\mathcal{C}}$ intersects  $\Gamma_{\theta(\phi)}$  in the two points contained in the edges $e$ and $f$.    

Using the hypothesis applied to $c$ and applying Lemma \ref{annuluslemma}, we can conclude that $H_1\cup_{c}H_2$ is a handlebody.    To complete the construction of the complement of $\Gamma_{\theta(\phi)}$, we attach 1-handles to $H_1\cup_{c}H_2$ corresponding to the remaining disk regions of $\Gamma-{\theta(\phi)}$ in $F$.   Thus for any $e,f$ in $\Gamma$, $(\Gamma-(e\cup{f}))_{\theta(\phi)}$ is free, as required.

\section{Closing question}
\label{question}

Since all graphs have freeness index at least 1, and since the freeness index of the Petersen graph is 4, we close with the following question:\\

\begin{Que}  
Do all graphs have freeness index at least 2? 3? 4?    
\end{Que}

We conjecture that every graph has freeness index at least 2, and that there exists some graph with freeness index equal to 2.     Note that if there exists a graph with freeness index one, it would provide (applying Theorem \ref{main}) a counter-example to the orientable cycle cover conjecture.

\bigskip

\noindent \address{Department of Mathematics \\ 
University of California, Davis\\
Davis, CA 95616}\\
\email{thompson@math.ucdavis.edu}

\end{document}